\documentclass[12pt,a4paper,onecolumn]{article}
\usepackage{textcomp}
\usepackage{tikz}
\usepackage{amssymb}
\usepackage{caption}
\usepackage{color}
\usepackage{fancybox}
\usepackage{verbatim}
\usepackage{fancybox}
\usepackage{titlesec} 
\usepackage{mathtools}
\usepackage{amsmath}
\usepackage{esvect}
\usepackage{cite}
\usepackage{setspace}
\usepackage{diagbox}
\numberwithin{equation}{section}
\usepackage[hidelinks]{hyperref}
\usepackage{mathrsfs}
\usetikzlibrary{positioning,fit,calc}
\usepackage{array}
\usepackage{bm}
\usepackage{multirow}
\usepackage{hhline}
\usepackage[super]{nth}
\usetikzlibrary{positioning,decorations.pathreplacing}
\usepackage{subcaption} 
\usetikzlibrary{decorations.pathmorphing}
\usetikzlibrary{arrows}
\usepackage{fixltx2e}
\usepackage{ragged2e}
\usepackage[super]{nth}

\usepackage{mathrsfs}
\usetikzlibrary{positioning,fit,calc}
\usepackage[a4paper,left=1in,top=1in,right=1in,bottom=1in,nohead]{geometry} 
\usepackage[hidelinks]{hyperref}
\usepackage{array}
\usepackage{bm}
\usepackage[symbol]{footmisc} 
\numberwithin{equation}{section}
\usepackage{amsthm}
\def\correspondingauthor{\footnote{Corresponding author. Email: williewong088@gmail.com.}}

\tikzset{block/.style={draw,thick,text width=2cm,minimum height=1cm,align=center},
         line/.style={-latex}}
\newcolumntype{P}[1]{>{\centering\arraybackslash}m{#1}} 

\titleformat{\section}[block]{\large\scshape\bfseries}{\thesection.}{1em}{} 
\titleformat{\subsection}[block]{\bfseries}{\thesubsection.}{1em}{} 

\DeclareFontFamily{U}{mathb}{}
\DeclareFontShape{U}{mathb}{m}{n}{<-5.5> mathb5 <5.5-6.5> mathb6 
<6.5-7.5> mathb7 <7.5-8.5> mathb8 <8.5-9.5> mathb9 <9.5-11> mathb10 
<11-> mathb12}{}
\DeclareSymbolFont{mathb}{U}{mathb}{m}{n}
\DeclareFontSubstitution{U}{mathb}{m}{n}

\DeclareMathSymbol{\blackdiamond}{\mathbin}{mathb}{"0C}

\newtheorem{thm}{Theorem}[section]
\newtheorem{ppn}[thm]{Proposition}
\newtheorem{cor}[thm]{Corollary}
\newtheorem{lem}[thm]{Lemma}

\theoremstyle{definition}
\newtheorem{defn}[thm]{Definition}

\newtheorem{rmk}[thm]{Remark}
\newtheorem{eg}[thm]{Example}

\newtheorem{customprob}{Problem}

\begin{document}
\pagenumbering{arabic}
\begin{center}
    \textbf{\Large Tournament completions of bipartite\\tournaments and their augmented directed cycles}
\vspace{0.1 in} 
    \\{\large H.W. Willie Wong\correspondingauthor{}}
\vspace{0.1 in} 
\\National Institute of Education\\Nanyang Technological University\\Singapore
\end{center}

\begin{abstract}
\noindent A tournament $T$ is a tournament completion of a bipartite tournament $D$ if $D$ is a spanning subdigraph of $T$, i.e., $V(D)=V(T)$ and $A(D)\subseteq A(T)$. If $C$ is a $k$-dicycle (i.e., directed cycle of length $k$) in a tournament completion $T$ of $D$ and $C$ is not a dicycle in $D$, i.e., $A(C)\subseteq A(T)$ and $A(C)\not\subseteq A(D)$, then we call $C$ an augmented $k$-dicycle of $T$. In this paper, we investigate the families of bipartite tournaments for which there exists a tournament completion with exactly one augmented $3$-dicycle and with no augmented $4$-dicycles. Our investigation may be viewed as a variant of the orientation completion problem initiated by Bang-Jensen et al..
\end{abstract}
\textbf{Keywords}: bipartite tournament; tournament completion; augmented directed cycle
\section{Introduction}
\begin{figure}[h]
\begin{subfigure}{.6\textwidth}
\begin{center}
\tikzstyle{every node}=[circle, draw, fill=black!100,
                       inner sep=0pt, minimum width=5pt]
\begin{tikzpicture}[thick,scale=0.6]%
\draw(2,0)node[label={[yshift=0cm]270:{$u_1$}}](u1){};
\draw(4,0)node[label={[yshift=0cm]270:{$u_2$}}](u2){};
\draw(0,2)node[label={[yshift=0cm]90:{$v_1$}}](v1){};
\draw(2,2)node[label={[yshift=0cm]90:{$v_2$}}](v2){};
\draw(4,2)node[label={[yshift=0cm]90:{$v_3$}}](v3){};
\draw(6,2)node[label={[yshift=0cm]90:{$v_4$}}](v4){};

\draw[->, line width=0.3mm, >=latex, shorten <= 0.1cm, shorten >= 0.1cm](u1)--(v1);
\draw[->, line width=0.3mm, >=latex, shorten <= 0.1cm, shorten >= 0.1cm](u1)--(v2);
\draw[->, line width=0.3mm, >=latex, shorten <= 0.1cm, shorten >= 0.1cm](v1)--(u2);
\draw[->, line width=0.3mm, >=latex, shorten <= 0.1cm, shorten >= 0.1cm](v2)--(u2);

\draw[->, line width=0.3mm, >=latex, shorten <= 0.1cm, shorten >= 0.1cm](u2)--(v3);
\draw[->, line width=0.3mm, >=latex, shorten <= 0.1cm, shorten >= 0.1cm](u2)--(v4);
\draw[->, line width=0.3mm, >=latex, shorten <= 0.1cm, shorten >= 0.1cm](v3)--(u1);
\draw[->, line width=0.3mm, >=latex, shorten <= 0.1cm, shorten >= 0.1cm](v4)--(u1);
\end{tikzpicture}
{\caption{$F$}}
\end{center}
\end{subfigure}%
\begin{subfigure}{.3\textwidth}
\begin{center}
\tikzstyle{every node}=[circle, draw, fill=black!100,
                       inner sep=0pt, minimum width=5pt]
\begin{tikzpicture}[thick,scale=0.6]%
\draw(2,0)node[label={[yshift=0cm]270:{$u_1$}}](u1){};
\draw(4,0)node[label={[yshift=0cm]270:{$u_2$}}](u2){};
\draw(0,2)node[label={[xshift=-0.1cm]90:{$v_1$}}](v1){};
\draw(2,2)node[label={[yshift=0cm]90:{$v_2$}}](v2){};
\draw(4,2)node[label={[yshift=0cm]90:{$v_3$}}](v3){};
\draw(6,2)node[label={[xshift=0.1cm]90:{$v_4$}}](v4){};

\draw[->, line width=0.3mm, >=latex, shorten <= 0.1cm, shorten >= 0.1cm](u1)--(v1);
\draw[->, line width=0.3mm, >=latex, shorten <= 0.1cm, shorten >= 0.1cm](u1)--(v2);
\draw[->, line width=0.3mm, >=latex, shorten <= 0.1cm, shorten >= 0.1cm](v1)--(u2);
\draw[->, line width=0.3mm, >=latex, shorten <= 0.1cm, shorten >= 0.1cm](v2)--(u2);

\draw[->, line width=0.3mm, >=latex, shorten <= 0.1cm, shorten >= 0.1cm](u2)--(v3);
\draw[->, line width=0.3mm, >=latex, shorten <= 0.1cm, shorten >= 0.1cm](u2)--(v4);
\draw[->, line width=0.3mm, >=latex, shorten <= 0.1cm, shorten >= 0.1cm](v3)--(u1);
\draw[->, line width=0.3mm, >=latex, shorten <= 0.1cm, shorten >= 0.1cm](v4)--(u1);

\draw[dashed, ->, line width=0.3mm, >=latex, shorten <= 0.1cm, shorten >= 0.1cm](u1)--(u2);

\draw[dashed, ->, line width=0.3mm, >=latex, shorten <= 0.1cm, shorten >= 0.1cm](v4)--(v3);
\draw[dashed, ->, line width=0.3mm, >=latex, shorten <= 0.1cm, shorten >= 0.1cm](v2)--(v3);
\draw[dashed, ->, line width=0.3mm, >=latex, shorten <= 0.1cm, shorten >= 0.1cm](v1)--(v2);

\draw[dashed, ->, line width=0.3mm, >=latex, shorten <= 0.1cm, shorten >= 0.1cm](v1) to [out=45, in=135](v3);
\draw[dashed, ->, line width=0.3mm, >=latex, shorten <= 0.1cm, shorten >= 0.1cm](v2) to [out=45, in=135](v4);
\draw[dashed, ->, line width=0.3mm, >=latex, shorten <= 0.1cm, shorten >= 0.1cm](v4) to [out=120, in=60](v1);
\end{tikzpicture}
{\caption{$T$}}
\end{center}
\end{subfigure}\hfill%
{\caption{A tournament completion $T$ of a bipartite tournament $F$.}\label{figA1.1}}
\end{figure}
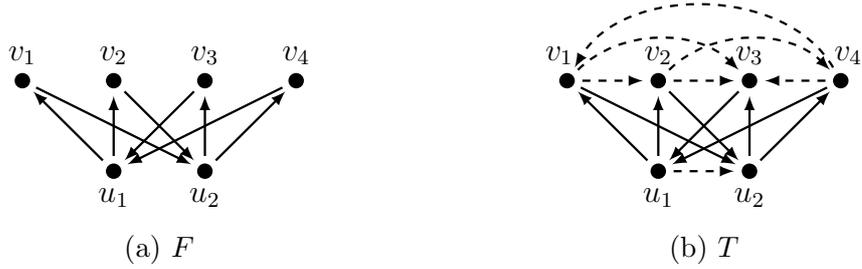

Let $G$ be a graph with vertex set $V(G)$ and edge set $E(G)$. For a digraph $D$, we denote its vertex set by $V(D)$ and its arc set by $A(D)$. We say that $D$ \textit{contains} a digraph $D'$ if $D'$ is a subdigraph of $D$, i.e., $V(D')\subseteq V(D)$ and $A(D')\subseteq A(D)$. The subdigraph of $D$ induced by the set of vertices $V\subseteq V(D)$ (resp. set of arcs $A\subseteq A(D)$) is denoted by $D\langle V\rangle$ (resp. $D\langle A\rangle$). An $\textit{orientation}$ of a graph $G$ is a digraph obtained from $G$ by assigning a direction to every edge in $G$. A \textit{tournament} is an orientation of a complete graph. A \textit{bipartite tournament} is an orientation of a complete bipartite graph; and we refer to its two nonempty partite sets as $V_1$ and $V_2$. A tournament $T$ is a \textit{tournament completion} of a bipartite tournament $D$ if $D$ is a spanning subdigraph of $T$, i.e., $V(D)=V(T)$ and $A(D)\subseteq A(T)$. If $C$ is a $k$-dicycle (i.e., directed cycle of length $k$) in a tournament completion $T$ of $D$ and $C$ is not a dicycle in $D$, i.e., $A(C)\subseteq A(T)$ and $A(C)\not\subseteq A(D)$, then we call $C$ an \textit{augmented $k$-dicycle} of $T$. If $\mathcal{K}\subseteq \{(j,k-j)\mid j\in\mathbb{Z}^+, \lceil\frac{k}{2}\rceil\le j\le k\}$, $|V(C)\cap V_i|=j$, and $|V(C)\cap V_{3-i}|=k-j$ for some $(j,k-j)\in \mathcal{K}$ and some $i=1,2$, then we call $C$ an \textit{augmented $\mathcal{K}_D$-dicycle}. In particular, if $\mathcal{K}=\{(j^*, k-j^*)\}$, then we may refer to an augmented $\mathcal{K}_D$-dicycle $C$ simply as an augmented $(j^*,k-j^*)_D$-dicycle. Whenever there is no ambiguity, we omit the subscript $D$ for the above notations. We refer the reader to \cite{JB GG} for any undefined terminology.
\begin{eg} Let $F$ be the bipartite tournament and $T$ be the tournament completion of $F$ shown in Figures \ref{figA1.1}(a) and \ref{figA1.1}(b), respectively. Then, $v_1 v_2 v_4 v_1$ is an augmented $(3,0)_F$-dicycle; $u_1 u_2 v_3 u_1$ and $v_2 v_3 u_1 v_2$ are augmented $(2,1)_F$-dicycles; and $u_1 u_2 v_4 v_3 u_1$ is an augmented $(2,2)_F$-dicycle. The $4$-dicycle $u_1 v_2 u_2 v_3 u_1$ is not an augmented $4$-dicycle since all its arcs are in $A(F)$.
\end{eg}

\noindent\par Tournaments and bipartite tournaments are among the most well-studied classes of digraphs. In particular, the topic of directed cycles in these two classes of digraphs has been extensively researched on (see \cite{JWM2, FH LM, KSB LWB, BA CT, NL AM, LWB CHCL, LV}). For example, it is well-known that a tournament (resp. bipartite tournament) contains a dicycle if and only if it contains a $3$-dicycle (resp. $4$-dicycle). However, to the author's knowledge, there has not been much research on tournament completions of bipartite tournaments and their augmented dicycles; and this paper aims to fill this research gap. Hence, we propose the following problem.
\begin{customprob} \label{probA}
Let $t\ge 0$ and $k\ge 3$ be given integers and let $\mathcal{K}\subseteq \{(j,k-j)\mid j\in\mathbb{Z}^+, \lceil\frac{k}{2}\rceil\le j\le k\}$. Characterize the family $\mathcal{F}$ of bipartite tournaments $D$ such that for any $D\in \mathcal{F}$, there exists some tournament completion of $D$ with exactly $t$ augmented $\mathcal{K}$-dicycles.
\end{customprob}

\begin{eg} \label{egA1.2}
If $t=0$, $k=3$ and $\mathcal{K}=\{(3,0), (2,1)\}$, then Problem \ref{probA} asks for the family $\mathcal{F}$ of bipartite tournaments $D$ for which there exists a tournament completion $T$ of $D$ with no augmented $3$-dicycles. Recall that a bipartite tournament contains no odd dicycles, particularly, no $3$-dicycles; and that a tournament is acyclic if and only if it contains no $3$-dicycles. Hence, any $3$-dicycle in $T$ is an augmented $3$-dicycle. Since a partially oriented graph can be completed to an acyclic-oriented graph if and only if it is acyclic (see \cite{JB JH XZ}), it follows that $\mathcal{F}$ is precisely the family of acyclic bipartite tournaments.
\end{eg}

\begin{eg}\label{egA1.3}
For any bipartite tournament $D$, there exists a tournament completion of $D$ with no augmented $(k,0)$-dicycles for all $k\ge 3$. This corresponds to the case of Problem \ref{probA} where $t=0$, $k\ge 3$, and $\mathcal{K}=\{(k,0)\}$. We sketch the proof as follows. For each $i=1,2$, linearly order the vertices in $V_i$ and define a tournament completion $T$ of $D$ such that for any $u,v\in V_i$, $uv\in A(T)$ if and only if $u$ precedes $v$. It is straightforward to see that both $T\langle V_1\rangle$ and $T\langle V_2\rangle$ are transitive tournaments and hence contain no augmented $(k,0)$-dicycles.
\end{eg}

\noindent\par In this paper, we determine the family $\mathcal{F}$ in Problem \ref{probA} for cases with a focus on augmented $3$-dicycles and $4$-dicycles. We note that the case in which $t=0$, $k=3$ and $\mathcal{K}=\{(2,1)\}$ will be resolved later as Proposition \ref{ppnA3.2}. Hence, in view of Examples \ref{egA1.2} and \ref{egA1.3}, we consider $t=1$ for augmented $3$-dicycles and $t=0$ for augmented $4$-dicycles. We summarize our results in Table \ref{tabA1.1}.
~\\~\\
\begin{spacing}{1.1}
\centering
\begin{tabular}{| P{2cm} | P{2cm} | P{5 cm}| P{4 cm}|}
\hline
$\bm{t}$ & $\bm{k}$ & $\bm{\mathcal{K}}$ & \textbf{Reference}\\
\hline
$1$ & $3$ & $\{(2,1)\}$ & Theorem \ref{thmA3.4} \\
\hline
$1$ & $3$ & $\{(3,0), (2,1)\}$ & Theorem \ref{thmA3.7} \\ 
\hline
$0$ & $4$ & $\{(2,2)\}$ & Proposition \ref{ppnA4.7} \\
\hline
$0$ & $4$ & $\{(3,1)\}$ & Proposition \ref{ppnA4.11} \\
\hline
$0$ & $4$ & $\{(3,1), (2,2)\}$ & Theorem \ref{thmA4.13} \\
\hline
\end{tabular}
{\captionof{table}{Summary of cases of Problem \ref{probA}.}\label{tabA1.1}}
\end{spacing}

\noindent\par The orientation completion problem was initiated by Bang-Jensen et al. \cite{JB JH XZ} and asked for a fixed class $\mathscr{C}$ of oriented graphs whether a given partially oriented graph can be completed to an oriented graph belonging to the class by orienting the (nonoriented) edges. The problem commonly generalizes orientation and recognition problems for graph and digraph classes and extends the representation extension problem. Bang-Jensen et al. \cite{JB JH XZ} and Hsu and Huang \cite{KH and JH} examined the problem for classes such as the local tournaments, locally transitive tournaments, and acyclic local tournaments. If $k\ge 3$ is odd and $\mathcal{K}=\{(j,k-j)\mid j\in\mathbb{Z}^+, \lceil\frac{k}{2}\rceil\le j\le k\}$, then since a bipartite tournament contains no odd dicycles, Problem \ref{probA} is simply the orientation completion problem asking for the fixed class $\mathscr{C}$ of tournaments (completions) with exactly $t$ $k$-dicycles and the given partially oriented graph being a bipartite tournament. However, if $k\ge 4$ is even and $\mathcal{K}=\{(j,k-j)\mid j\in\mathbb{Z}^+, \lceil\frac{k}{2}\rceil\le j\le k\}$, the class of tournament completions with exactly $t$ augmented $k$-dicycles is not fixed as they may contain a different number of $k$-dicycles, which depends on their given bipartite tournaments.  Hence, Problem \ref{probA} may be viewed as a variant of the orientation completion problem.

\noindent\par This paper is organised as follows. In Section 2, we provide the necessary notation and terminology, and some useful results for Section 3. In Sections 3 and 4, we prove our results on augmented $3$-dicycles and $4$-dicycles, respectively. In the last section, we give some concluding remarks.
 
\section{Preliminaries}
In this paper, let $[n]=\{1,2,\ldots, n\}$ for any $n\in\mathbb{Z}^+$, and $[n]=\emptyset$ otherwise. Let $D$ be a digraph.  If $uv\in A(D)$ for some $u,v\in V(D)$, then we denote it by $u\rightarrow v$. For a vertex $v\in V(D)$, the \textit{outset} and \textit{inset} of $v$ are defined to be $N^+_D(v)=\{x\in V(D)\mid v\rightarrow x\}$ and $N^-_D(v)=\{y\in V(D)\mid y\rightarrow v\}$, respectively. The \textit{outdegree} and \textit{indegree} of $v$ are defined by $\deg^+_D(v)=|N^+_D(v)|$ and $\deg^-_D(v)=|N^-_D(v)|$, respectively. If $V\subseteq V(D)$ and $V\subseteq N^+(v)$ (resp. $V\subseteq N^-(v)$), then we denote this by $v\rightarrow V$ (resp. $V\rightarrow v$). For any vertices $u,v\in V(D)$, the $\textit{distance}$ $d_D(u,v)$ from $u$ to $v$ is defined as the length of a shortest directed path (i.e., number of arcs) from $u$ to $v$. If there is no ambiguity, we shall omit the subscript $D$ for the above notations.
\noindent\par We will make frequent use of a characterization of acyclic bipartite tournaments, namely Theorem \ref{thmA2.2}. Given a finite nonempty set $X\subseteq \mathbb{Z}^+$, let $D_X$ be the digraph with 
\begin{align}
V(D_X)=X \text{ and }A(D_X)=\{ab\mid a<b \text{ and } a\not\equiv b\pmod{2}\}. \label{eqA2.1}
\end{align}
Clearly, $D_X$ is a bipartite tournament.
\begin{lem}\label{lemA2.1}
Let $X=\{x_1, x_2 \ldots, x_n\}\subset \mathbb{Z}^+$ for some positive integer $n$. If $x_i$ and $x_j$ are in the same partite set and $\deg^+_{D_X}(x_i)>\deg^+_{D_X}(x_j)$ for some $x_i, x_j\in X$, then $x_i<x_j$. 
\end{lem}
\noindent\textit{Proof}: If $x_i\ge x_j$, then $N^+_{D_X}(x_i)\subseteq N^+_{D_X}(x_j)$ by (\ref{eqA2.1}), which contradicts the fact that $\deg^+_{D_X}(x_i)>\deg^+_{D_X}(x_j)$.
\qed

\indent\par Das et al. \cite{SD PG SG SS} introduced the notion of $D_X$ and said that an oriented bipartite graph $D$ is \textit{bitransitive} if $\{v_1v_2, v_2v_3, v_3v_4\}\subseteq A(D)$ implies $v_1v_4\in A(D)$ for any $v_1, v_2, v_3, v_4\in V(D)$. With these notions, they proved the following.
\begin{thm} (Das et al. \cite{SD PG SG SS}) \label{thmA2.2} Let $D$ be a bipartite tournament. Then, the following are equivalent.
\\(1) $D$ is acyclic.
\\(2) $D$ contains no $4$-dicycle.
\\(3) $D\cong D_X$ for some nonempty set $X\subseteq \mathbb{Z}^+$.
\\(4) $D$ is bitransitive.
\end{thm}

\section{Exactly one augmented $3$-dicycle}
In this section, we characterize the bipartite tournaments for which there exists some tournament completion with exactly one augmented $3$-dicycle. In particular, we consider $\mathcal{K}=\{(2,1)\}$ and $\mathcal{K}=\{(3,0), (2,1)\}$ (and $t=1$ and $k=3$ as denoted in Problem \ref{probA}) in Theorems \ref{thmA3.4} and \ref{thmA3.7}, respectively.

\begin{lem}\label{lemA3.1}
Let $D$ be a bipartite tournament. If $D$ contains a dicycle, then every tournament completion of $D$ contains at least two augmented $(2,1)$-dicycles.
\end{lem}
\noindent\textit{Proof}: Since $D$ contains a dicycle, $D$ contains a $4$-dicycle, say $u_1 u_2 u_3 u_4 u_1$ by Theorem \ref{thmA2.2}. Any tournament completion $T$ of $D$ satisfies either $u_1u_3\in A(T)$ or $u_3u_1\in A(T)$. If the former (resp. latter) holds, then $u_1u_3 u_4 u_1$ (resp. $u_3u_1u_2u_3$) is an augmented $(2,1)$-dicycle in $T$. Similar case consideration can be applied to the arc $u_2u_4$ or $u_4u_2$ to prove the existence of another augmented $(2,1)$-dicycle in $T$.
\qed

\begin{ppn}\label{ppnA3.2}
Let $D$ be a bipartite tournament. There exists a tournament completion of $D$ with no augmented $(2,1)$-dicycles if and only if $D$ is acyclic.
\end{ppn}
\noindent\textit{Proof}: $(\Rightarrow)$ This follows directly from Lemma \ref{lemA3.1}.
\noindent\par $(\Leftarrow)$ By Example \ref{egA1.2}, there exists a tournament completion of $D$ with no augmented $3$-dicycles, in particular, no augmented $(2,1)$-dicycles.
\qed

\indent\par Before we prove Theorem \ref{thmA3.4}, we characterize the tournaments of order $n$ with exactly one dicycle. For any given $n\ge 3$ and $1\le r \le n-2$, let $T^r_n$ be the tournament where $V(T^r_n)=[n]$ and $ij\in A(T^r_n)$ if and only if $i<j$ and $(i,j)\neq (r,r+2)$. In particular, we have $[r-1]\rightarrow [n]\backslash [r-1]$ and $\{r, r+1, r+2\}\rightarrow [n]\backslash [r+2]$ and $r\rightarrow r+1\rightarrow r+2\rightarrow r$ (see Figure \ref{figA3.2} for $T^3_8$).
\begin{center}
\tikzstyle{every node}=[circle, draw, fill=black!100,
                       inner sep=0pt, minimum width=6pt]
\begin{tikzpicture}[thick,scale=0.7]%
\draw(0,1)node[label={[yshift=0.1cm] 90:{$3$}}](3){};
\draw(-0.866,-0.5)node[label={[yshift=-0.08cm, xshift=-0.1cm] 225:{$4$}}](4){};
\draw(0.866,-0.5)node[label={[yshift=-0.08cm, xshift=0.1cm] 315:{$5$}}](5){};

\draw[->, line width=0.3mm, >=latex, shorten <= 0.1cm, shorten >= 0.1cm](3)--(4);
\draw[->, line width=0.3mm, >=latex, shorten <= 0.1cm, shorten >= 0.1cm](4)--(5);
\draw[->, line width=0.3mm, >=latex, shorten <= 0.1cm, shorten >= 0.1cm](5)--(3);

\draw(-6.5,0.2)node[label={[yshift=-0.3, xshift=-0.1cm]150:{$1$}}](1){};
\draw(-5.5,-0.2)node[label={[yshift=0.35, xshift=0.1cm]330:{$2$}}](2){};
\draw[->, line width=0.3mm, >=latex, shorten <= 0.05cm, shorten >= 0.01cm](1)--(2);

\draw (-6,0) node[circle, draw, fill=none, inner sep=0pt, minimum width=60pt, label=270:{$[2]$}]([2]){};

\draw(6.65,0.1)node[label={[yshift=0cm, xshift=0.05cm]0:{$8$}}](8){};
\draw(5.65,-0.5)node[label={[yshift=0cm, xshift=-0.05cm]210:{$7$}}](7){};
\draw(5.65,0.5)node[label={[yshift=0cm, xshift=-0.05cm]150:{$6$}}](6){};
\draw[->, line width=0.3mm, >=latex, shorten <= 0.05cm, shorten >= 0.01cm](6)--(7);
\draw[->, line width=0.3mm, >=latex, shorten <= 0.05cm, shorten >= 0.01cm](7)--(8);
\draw[->, line width=0.3mm, >=latex, shorten <= 0.05cm, shorten >= 0.01cm](6)--(8);
\draw (6,0) node[circle, draw, fill=none, inner sep=0pt, minimum width=60pt, label={[yshift=0.3cm] 270:{$[8]\backslash [5]$}}](n){};

\draw[line,shorten <= 0.65cm, shorten >= 0.65cm] (-5,0) to (-1,0);
\draw[line,shorten <= 0.65cm, shorten >= 0.65cm] (1,0) to (5,0);
\draw[line,shorten <= 0.65cm, shorten >= 0.65cm] (-4.75,0.5) to [out=35, in=145] (4.75,0.5);
\end{tikzpicture}
{\captionof{figure}{$T^3_8$.}\label{figA3.2}}
\end{center}
\begin{lem}\label{lemA3.3}
A tournament $T$ of order $n\ge 3$ has exactly one dicycle if and only if $T\cong T^r_n$ for some $r\in [n-2]$.
\end{lem}
\noindent\textit{Proof}: $(\Leftarrow)$ This follows directly from the definition of $T^r_n$.
\noindent\par $(\Rightarrow)$ If $n=3$, then $T\cong T^1_3$ is simply the dicycle of length $3$. Now, let $T$ be a tournament of order $n\ge 4$ with exactly one dicycle $C$. Since a tournament contains a dicycle if and only if it contains a $3$-dicycle, $C$ must be a dicycle of length $3$, say $C=u_1 u_2 u_3 u_1$. So, $T-u_1$ is a transitive tournament, and there exists a bijection $\sigma:V(T-u_1)\rightarrow [n-1]$ such that $u_i u_j\in A(T-u_1)$ if and only if $\sigma(u_i)<\sigma(u_j)$. Without loss of generality, we may assume $\sigma(u_2)=r$ for some integer $r\in [n-1]$. Then, $\sigma(u_3)>r$ since $u_2\rightarrow u_3$. Furthermore, $\sigma(u_3)=r+1$; for if $\sigma(u_3)=r+i$ for some $i\ge 2$, then $T$ contains another dicycle, namely $u_2 v u_3 u_1 u_2$, where $v$ is a vertex with $\sigma(v)=r+1$, a contradiction. Hence, $\{\sigma(u)\mid u\in N^+_T(u_2)\backslash V(C)\}=\{\sigma(u)\mid u\in N^+_T(u_3)\backslash V(C)\}=[n-1]\backslash [r+1]$ and $\{\sigma(u)\mid N^-_T(u_2)\backslash V(C)\}=\{\sigma(u)\mid N^-_T(u_3)\backslash V(C)\}=[r-1]$.
\noindent\par We claim that $N^+_T(u_1)\backslash V(C)=N^+_T(u_2)\backslash V(C)$, or equivalently, $N^-_T(u_1)\backslash V(C)=N^-_T(u_2)\backslash V(C)$. If $v\in (N^+_T(u_1)\backslash V(C))\cap (N^-_T(u_2)\backslash V(C))$, then $u_1 v u_2 u_3 u_1$ is another dicycle in $T$, a contradiction. If $v\in (N^+_T(u_2)\backslash V(C))\cap (N^-_T(u_1)\backslash V(C))$, then $u_2 v u_1 u_2$ is another dicycle in $T$, a contradiction. Hence, the claim follows.
\noindent\par Now, let $\sigma^*:V(T)\rightarrow [n]$, where  
\begin{align*}
\sigma^*(u)= \left\{
  \begin{array}{@{}ll@{}}
    r,& \text{if $u=u_1$}, \nonumber\\
    \sigma(u)+1,& \text{if $u\in V(T-u_1)$ and $\sigma(u)\ge r$}, \nonumber\\
    \sigma(u),& \text{if $u\in V(T-u_1)$ and $\sigma(u)< r$}. \nonumber
  \end{array}\right.
\end{align*}
It is straightforward to see from the bijection $\sigma^*$ that $T\cong T^r_{n}$.
\qed

\begin{thm}\label{thmA3.4}
Let $D$ be a bipartite tournament of order $n$. Then, the following are equivalent.
\\(1) There exists a tournament completion of $D$ with exactly one augmented $(2,1)$-dicycle.
\\(2) $D$ is acyclic and there exist two vertices $u_1$ and $u_2$ belonging to a common partite set such that $\deg^+_D(u_1)=\deg^+_D(u_2)+1$.
\\(3) There exists some $X=\{x_1, x_2, \ldots, x_n\}\subset \mathbb{Z}^+$ such that $D\cong D_X$, $x_1<x_2<\cdots<x_n$ and $x_l\equiv x_{l+2}\not\equiv x_{l+1}\pmod{2}$ for some $l\in [n-2]$.
\end{thm}
\noindent\textit{Proof}: ((1)$\Rightarrow$(2)) By Lemma \ref{lemA3.1}, $D$ is acyclic. Let $T$ be a tournament completion of $D$ with exactly one augmented $(2,1)$-dicycle, say $u_1 v u_2 u_1$, where $u_1, u_2\in V_i$ and $v\in V_{3-i}$ for exactly one of $i=1,2$. Now, if there exists some vertex $w\in (N^+_D(u_1)\cap N^-_D(u_2))\backslash\{v\}$, then $u_1 w u_2 u_1$ is another augmented $(2,1)$-dicycle in $T$, a contradiction. Hence, $N^+_D(u_1)\cap N^-_D(u_2)=\{v\}$. By Theorem \ref{thmA2.2} and the acyclicity of $D$, $D$ is bitransitive. Since $N^+_D(u_2)\subseteq N^+_D(u_1)$ by the bitransitivity property of $D$, it follows that $\deg^+_D(u_1)=\deg^+_D(u_2)+1$.
\begin{figure}[h]
\begin{center}
\tikzstyle{every node}=[circle, draw, fill=black!100,
                       inner sep=0pt, minimum width=5pt]
\begin{tikzpicture}[thick,scale=0.6]%
\draw(3,0)node[label={[yshift=-0.1cm]270:{$x_l$}}](l){};
\draw(0,0)node[label={[yshift=-0.1cm]270:{$x_i$}}](i){};
\draw(1,0)node[label={[yshift=-0.1cm]270:{$$}}](i_1){};

\draw(8,5)node[label={[yshift=-0.15cm]90:{$x_{s-1}$}}](s-1){};
\draw(5,5)node[label={[yshift=-0.15cm]90:{$x_{l+1}$}}](l+1){};
\draw(6,5)node[label={[yshift=-0.15cm]90:{$$}}](l+1_1){};

\draw(13,0)node[label={[yshift=-0.1cm]270:{$x_j$}}](j){};
\draw(10,0)node[label={[yshift=-0.1cm]270:{$x_s$}}](s){};
\draw(11,0)node[label={[yshift=-0.1cm]270:{$$}}](s_1){};

\draw(-3,5)node[label={[yshift=-0.1cm]270:{$$}}](A1){};
\draw(-2,5)node[label={[yshift=-0.1cm]270:{$$}}](A2){};

\draw(15,5)node[label={[yshift=-0.1cm]270:{$$}}](A3){};
\draw(16,5)node[label={[yshift=-0.1cm]270:{$$}}](A4){};

\node[rectangle, draw, fill=none, inner xsep=3mm,inner ysep=5mm, fit=(i)(l), label={[xshift=0cm, yshift=0.75cm] 270:$\{x_a\mid a\in\mathscr{A}\}$}](R1){};
\node[rectangle, draw, fill=none, inner xsep=3mm,inner ysep=5mm, fit=(s)(j), label={[xshift=0cm, yshift=0.75cm] 270:$\{x_b\mid b\in\mathscr{B}\}$}](R3){};
\node[rectangle, draw, fill=none, inner xsep=3mm,inner ysep=5mm, fit=(l+1)(s-1), label={[xshift=0cm, yshift=0.75cm] 270:$$}](R2){};

\draw[->, line width=0.3mm, >=latex, shorten <= 0.2cm, shorten >= 0.2cm](R1.north)--(R2);
\draw[dotted, ->, line width=0.3mm, >=latex, shorten <= 0.2cm, shorten >= 0.2cm](R2)--(R3.north);

\draw[dotted, ->, line width=0.3mm, >=latex, shorten <= 0.2cm, shorten >= 0.2cm](A1)--(R1);
\draw[dotted, ->, line width=0.3mm, >=latex, shorten <= 0.2cm, shorten >= 0.2cm](A2)--(R1);
\draw[dotted, ->, line width=0.3mm, >=latex, shorten <= 0.2cm, shorten >= 0.2cm](A1)--(R3.north);
\draw[dotted, ->, line width=0.3mm, >=latex, shorten <= 0.2cm, shorten >= 0.2cm](A2)--(R3.north);

\draw[->, line width=0.3mm, >=latex, shorten <= 0.2cm, shorten >= 0.2cm](R1.north)--(A3);
\draw[->, line width=0.3mm, >=latex, shorten <= 0.2cm, shorten >= 0.2cm](R1.north)--(A4);
\draw[->, line width=0.3mm, >=latex, shorten <= 0.2cm, shorten >= 0.2cm](R3)--(A3);
\draw[->, line width=0.3mm, >=latex, shorten <= 0.2cm, shorten >= 0.2cm](R3)--(A4);

\draw[loosely dotted, line width=0.3mm, >=latex, shorten <= 0.2cm, shorten >= 0.15cm](i_1)--(l);
\draw[loosely dotted, line width=0.3mm, >=latex, shorten <= 0.2cm, shorten >= 0.15cm](l+1_1)--(s-1);
\draw[loosely dotted, line width=0.3mm, >=latex, shorten <= 0.2cm, shorten >= 0.15cm](s_1)--(j);
\end{tikzpicture}
{\caption{An illustration of $x_l$ and $x_s$ for the proof of ((2)$\Rightarrow$(3)) in Theorem \ref{thmA3.4}.}\label{figA3.3}}
\end{center}
\end{figure}
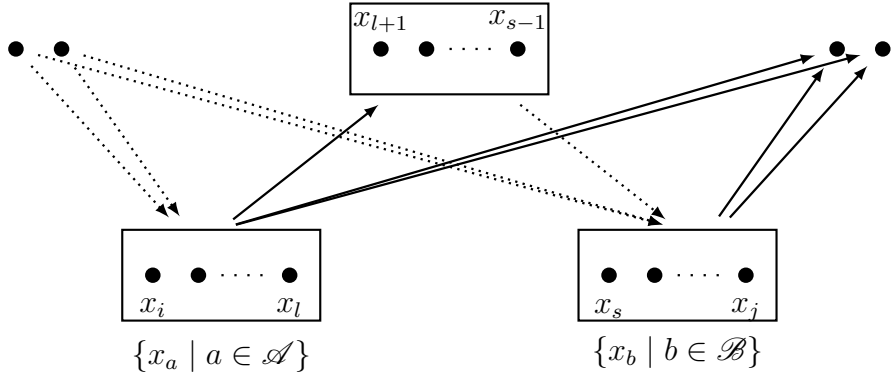

\noindent\par ((2)$\Rightarrow$(3)) Since $D$ is acyclic, by Theorem \ref{thmA2.2}, there exists some $X=\{x_1, x_2, \ldots, x_n\}$ $\subset \mathbb{Z}^+$ such that $D\cong D_X$ and $x_1<x_2<\cdots<x_n$. Note that for any $i,j\in [n]$, $x_i<x_j$ if and only if $i<j$; we shall apply this understanding to any use of the subscripts $i$, $j$ of $x$ for the rest of the proof. Suppose the vertex $u_1$ (resp. $u_2$) corresponds to $x_i\in X$ (resp. $x_j\in X$) for some $i\in [n]$ (resp. $j\in [n]$). By the fact that $\deg^+_{D_X}(x_i)=\deg^+_{D_X}(x_j)+1$ and Lemma \ref{lemA2.1}, we have $x_i<x_j$, or equivalently, $1\le i<j\le n$. Let 
\begin{align*}
\mathscr{A}=\{r\in [n]\backslash [i-1]\mid \deg^+_{D_X}(x_r)=\deg^+_{D_X}(x_i) \text{ and }x_r\equiv x_i\pmod{2}\}.
\end{align*}
Clearly, $i\in \mathscr{A}$. The vertices in the set $\{x_a\mid a\in \mathscr{A}\}$ share the same outneighbours with $x_i$ in $D_X$. Denote the largest element in $\mathscr{A}$ by $l$ so that $\deg^+_{D_X}(x_l)=\deg^+_{D_X}(x_i)$, $x_l\equiv x_i\pmod{2}$ and $1\le i\le l$ (see Figure \ref{figA3.3}). Note that it is possible that $i=l$. Now, it follows that
\begin{align}
x_l\not\equiv x_{l+1}\pmod{2}. \label{eqA3.1}
\end{align}
Otherwise, the definition (\ref{eqA2.1}) of $D_X$ and the fact that $\{x\in X\mid x_l<x<x_{l+1}\}=\emptyset$ imply that $\deg^+_{D_X}(x_{l+1})=\deg^+_{D_X}(x_l)=\deg^+_{D_X}(x_i)$, which contradicts the maximality of $l$. Note that it also follows from $\deg^+_{D_X}(x_l)=\deg^+_{D_X}(x_i)=\deg^+_{D_X}(x_j)+1$ and Lemma 2.1 that $l<j$.
\noindent\par Similarly, let
\begin{align*}
\mathscr{B}=\{r\in [j]\mid \deg^+_{D_X}(x_r)=\deg^+_{D_X}(x_j) \text{ and }x_r\equiv x_j\pmod{2}\}.
\end{align*}
Clearly, $j\in \mathscr{B}$. The vertices in the set $\{x_b\mid b\in \mathscr{B}\}$ share the same outneighbours with $x_j$ in $D_X$. Denote the smallest element in $\mathscr{B}$ by $s$ so that $\deg^+_{D_X}(x_s)=\deg^+_{D_X}(x_j)$, $x_s\equiv x_j\pmod{2}$ and $s\le j\le n$ (see Figure \ref{figA3.3}). Note that it is possible that $j=s$. Now, it follows that
\begin{align}
x_{s-1}\not\equiv x_s\pmod{2}. \label{eqA3.2}
\end{align}
Otherwise, the definition (\ref{eqA2.1}) of $D_X$ and the fact that $\{x\in X\mid x_{s-1}<x<x_s\}=\emptyset$ imply that $\deg^+_{D_X}(x_{s-1})=\deg^+_{D_X}(x_s)=\deg^+_{D_X}(x_j)$, a contradiction to the minimality of $s$. Since $u_1$ and $u_2$ are in the same partite set, $x_l\equiv x_i\equiv x_j\equiv x_s\pmod{2}$. Then, it follows from $\deg^+_{D_X}(x_s)+1=\deg^+_{D_X}(x_j)+1=\deg^+_{D_X}(x_i)=\deg^+_{D_X}(x_l)$ and Lemma 2.1 that $l<s$. Consequently, $1\le i\le l<s\le j\le n$.
\noindent\par  It remains to show that $s=l+2$. Since $\deg^+_{D_X}(x_l)=\deg^+_{D_X}(x_s)+1$, we have $|N^+_{D_X}(x_l)\cap N^-_{D_X}(x_s)|>0$ and $s\neq l+1$. It is necessary from the linear order of $X$ and the definition (\ref{eqA2.1}) of $D_X$ that $s\ge l+2$ and $N^+_{D_X}(x_s)\subseteq N^+_{D_X}(x_l)$. Now, if $s\ge l+3$, then (\ref{eqA3.1}) and (\ref{eqA3.2}) imply $\{x_{l+1}, x_{s-1}\}\subseteq N^+_{D_X}(x_l)\cap N^-_{D_X}(x_s)$ and $\deg^+_{D}(u_1)=\deg^+_{D_X}(x_l)\ge\deg^+_{D_X}(x_s)+2=\deg^+_{D}(u_2)+2$, a contradiction. Therefore, $s=l+2$ as desired.

\noindent\par ((3)$\Rightarrow$(1)) Let $T$ be the tournament completion of $D_X$ defined as follows: For any $x,y\in X$, \begin{align}
xy\in A(T) \iff x<y \text{ except when }x=x_l\text{ and } y=x_{l+2}. \label{eqA3.3}
\end{align}
By the definition (\ref{eqA2.1}) of $D_X$ and (\ref{eqA3.3}), any augmented $(2,1)$-dicycle in $T$ must contain the arc $x_{l+2} x_l$. Since $N^+_{D_X}(x_l)\cap N^-_{D_X}(x_{l+2})=\{x_{l+1}\}$ by the definition of $D_X$, $x_l x_{l+1} x_{l+2} x_l$ is the only augmented $(2,1)$-dicycle in $T$.
\qed

\indent\par
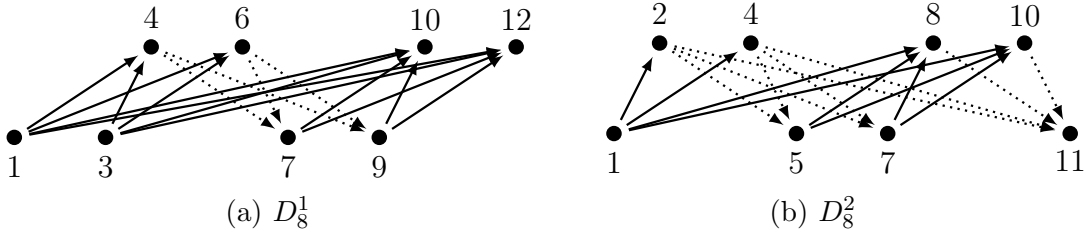
\begin{figure}[h]
\begin{subfigure}{.55\textwidth}
\begin{center}
\tikzstyle{every node}=[circle, draw, fill=black!100,
                       inner sep=0pt, minimum width=5pt]
\begin{tikzpicture}[thick,scale=0.6]%
\draw(8,0)node[label={[yshift=-0.1cm]270:{$9$}}](9){};
\draw(6,0)node[label={[yshift=-0.1cm]270:{$7$}}](7){};
\draw(2,0)node[label={[yshift=-0.1cm]270:{$3$}}](3){};
\draw(0,0)node[label={[yshift=-0.1cm]270:{$1$}}](1){};
\draw(11,2)node[label={[yshift=0cm]90:{$12$}}](12){};
\draw(9,2)node[label={[yshift=0cm]90:{$10$}}](10){};
\draw(5,2)node[label={[yshift=0.1cm]90:{$6$}}](6){};
\draw(3,2)node[label={[yshift=0.1cm]90:{$4$}}](4){};

\draw[->, line width=0.3mm, >=latex, shorten <= 0.1cm, shorten >= 0.1cm](1)--(4);
\draw[->, line width=0.3mm, >=latex, shorten <= 0.1cm, shorten >= 0.1cm](1)--(6);
\draw[->, line width=0.3mm, >=latex, shorten <= 0.1cm, shorten >= 0.1cm](3)--(4);
\draw[->, line width=0.3mm, >=latex, shorten <= 0.1cm, shorten >= 0.1cm](3)--(6);

\draw[dotted, ->, line width=0.3mm, >=latex, shorten <= 0.1cm, shorten >= 0.1cm](4)--(7);
\draw[dotted, ->, line width=0.3mm, >=latex, shorten <= 0.1cm, shorten >= 0.1cm](4)--(9);
\draw[dotted, ->, line width=0.3mm, >=latex, shorten <= 0.1cm, shorten >= 0.1cm](6)--(7);
\draw[dotted, ->, line width=0.3mm, >=latex, shorten <= 0.1cm, shorten >= 0.1cm](6)--(9);

\draw[->, line width=0.3mm, >=latex, shorten <= 0.1cm, shorten >= 0.1cm](7)--(10);
\draw[->, line width=0.3mm, >=latex, shorten <= 0.1cm, shorten >= 0.1cm](7)--(12);
\draw[->, line width=0.3mm, >=latex, shorten <= 0.1cm, shorten >= 0.1cm](9)--(10);
\draw[->, line width=0.3mm, >=latex, shorten <= 0.1cm, shorten >= 0.1cm](9)--(12);

\draw[->, line width=0.3mm, >=latex, shorten <= 0.1cm, shorten >= 0.1cm](1)--(10);
\draw[->, line width=0.3mm, >=latex, shorten <= 0.1cm, shorten >= 0.1cm](1)--(12);
\draw[->, line width=0.3mm, >=latex, shorten <= 0.1cm, shorten >= 0.1cm](3)--(10);
\draw[->, line width=0.3mm, >=latex, shorten <= 0.1cm, shorten >= 0.1cm](3)--(12);
\end{tikzpicture}
{\caption{$D^1_8$}}
\end{center}
\end{subfigure}%
\begin{subfigure}{.35\textwidth}
\begin{center}
\tikzstyle{every node}=[circle, draw, fill=black!100,
                       inner sep=0pt, minimum width=5pt]
\begin{tikzpicture}[thick,scale=0.6]%
\draw(10,0)node[label={[yshift=0cm]270:{$11$}}](11){};
\draw(6,0)node[label={[yshift=-0.1cm]270:{$7$}}](7){};
\draw(4,0)node[label={[yshift=-0.1cm]270:{$5$}}](5){};
\draw(0,0)node[label={[yshift=-0.1cm]270:{$1$}}](1){};
\draw(9,2)node[label={[yshift=0cm]90:{$10$}}](10){};
\draw(7,2)node[label={[yshift=0.1cm]90:{$8$}}](8){};
\draw(3,2)node[label={[yshift=0.1cm]90:{$4$}}](4){};
\draw(1,2)node[label={[yshift=0.1cm]90:{$2$}}](2){};

\draw[->, line width=0.3mm, >=latex, shorten <= 0.1cm, shorten >= 0.1cm](1)--(2);
\draw[->, line width=0.3mm, >=latex, shorten <= 0.1cm, shorten >= 0.1cm](1)--(4);
\draw[->, line width=0.3mm, >=latex, shorten <= 0.1cm, shorten >= 0.1cm](1)--(8);
\draw[->, line width=0.3mm, >=latex, shorten <= 0.1cm, shorten >= 0.1cm](1)--(10);

\draw[dotted, ->, line width=0.3mm, >=latex, shorten <= 0.1cm, shorten >= 0.1cm](2)--(5);
\draw[dotted, ->, line width=0.3mm, >=latex, shorten <= 0.1cm, shorten >= 0.1cm](4)--(5);
\draw[dotted, ->, line width=0.3mm, >=latex, shorten <= 0.1cm, shorten >= 0.1cm](2)--(7);
\draw[dotted, ->, line width=0.3mm, >=latex, shorten <= 0.1cm, shorten >= 0.1cm](4)--(7);
\draw[dotted, ->, line width=0.3mm, >=latex, shorten <= 0.1cm, shorten >= 0.1cm](2)--(11);
\draw[dotted, ->, line width=0.3mm, >=latex, shorten <= 0.1cm, shorten >= 0.1cm](4)--(11);

\draw[->, line width=0.3mm, >=latex, shorten <= 0.1cm, shorten >= 0.1cm](5)--(8);
\draw[->, line width=0.3mm, >=latex, shorten <= 0.1cm, shorten >= 0.1cm](5)--(10);
\draw[->, line width=0.3mm, >=latex, shorten <= 0.1cm, shorten >= 0.1cm](7)--(8);
\draw[->, line width=0.3mm, >=latex, shorten <= 0.1cm, shorten >= 0.1cm](7)--(10);

\draw[dotted, ->, line width=0.3mm, >=latex, shorten <= 0.1cm, shorten >= 0.1cm](8)--(11);
\draw[dotted, ->, line width=0.3mm, >=latex, shorten <= 0.1cm, shorten >= 0.1cm](10)--(11);
\end{tikzpicture}
{\caption{$D^2_8$}}
\end{center}
\end{subfigure}\hfill%
{\caption{Acyclic bipartite tournaments $D^1_8\cong D_{X^1_8}$ and $D^2_8\cong D_{X^2_8}$.}\label{figA3.4}}
\end{figure}

 For the case $t=1$, $k=3$, and $\mathcal{K}=\{(3,0), (2,1)\}$ of Problem \ref{probA}, we need to define two special forms of acyclic bipartite tournaments. For any given positive integer $n$, let $X^1_n$ (resp. $X^2_n$) be the set containing the $n$ smallest integers in $\{6i-5, 6i-3, 6i-2, 6i\mid i\in\mathbb{Z}^+\}$ (resp. $\{6i-5, 6i-4, 6i-2, 6i-1\mid i\in\mathbb{Z}^+\}$). For $j=1,2$, we define $D^j_n$ to be the acyclic bipartite tournament $D^j_n\cong D_{X^j_n}$ (see Figure \ref{figA3.4} when $n=8$).

\begin{lem}\label{lemA3.5}
Let $n\in \mathbb{Z}^+$ and $X=\{x_1, x_2, \ldots, x_n\}\subset \mathbb{Z}^+$, where $x_1<x_2<\cdots<x_n$. Suppose $D_X\cong D^i_n$ for some $i=1,2$. If $x\equiv y\equiv z\pmod{2}$ and $x<y<z$ for some $x,y,z\in X$, then $N^+_{D_X}(x)\cap N^-_{D_X}(z)\neq\emptyset$, and either $N^+_{D_X}(x)\cap N^-_{D_X}(y)\neq\emptyset$ or $N^+_{D_X}(y)\cap N^-_{D_X}(z)\neq\emptyset$.
\end{lem}
\noindent\textit{Proof}: We prove the case where $D_X\cong D^1_n$; the proof is similar if $D_X\cong D^2_n$. Suppose that $x\equiv y\equiv z\pmod{2}$ and $x<y<z$ for some $x,y,z\in X$. Observe that $x_j\equiv x_{j+1}\not\equiv x_{j+2}\pmod{2}$ for all odd $j\in [n-2]$. Consequently, there exists $a\in X$ such that $a\not\equiv x\pmod{2}$, and either $x<a<y<z$ or $x<y<a<z$. It follows that $a\in N^+_{D_X}(x)\cap N^-_{D_X}(z)$, and either $a\in N^+_{D_X}(x)\cap N^-_{D_X}(y)$ or $a\in N^+_{D_X}(y)\cap N^-_{D_X}(z)$, respectively.
\qed

\indent\par Lemma \ref{lemA3.6} provides an equivalent condition for all acyclic bipartite tournaments of order $n$ that are isomorphic to neither $D^1_n$ nor $D^2_n$. We note that an acyclic bipartite tournament is isomorphic to $D_X$ by Theorem \ref{thmA2.2}. 
\begin{lem}\label{lemA3.6}
Let $D_X$ be an acyclic bipartite tournament of order $n$, where $X=\{x_1, x_2,$ $\ldots, x_n\}\subset \mathbb{Z}^+$ and $x_1<x_2<\cdots<x_n$. Then, $D_X\not\cong D^i_n$ for all $i=1,2$, if and only if there exists some $j\in [n-2]$ such that $x_j\equiv x_{j+1}\equiv x_{j+2}\pmod{2}$ or $x_j\equiv x_{j+2}\not\equiv x_{j+1}\pmod{2}$.
\end{lem}
\noindent\textit{Proof}: $(\Leftarrow)$ By definition, each of $X^1_n$ and $X^2_n$ does not contain three consecutive elements of the same parity or of alternating parity.
\noindent\par $(\Rightarrow)$ We may assume that $x_1\equiv 1\pmod{2}$. For if $x_1\equiv 2\pmod{2}$, then we can consider $X'=\{x+1\mid x\in X\}$ in place of $X$, since $D_{X'}\cong D_X$. For each $r=1,2$, denote the $i$-th element of $X^r_n$ as $y^{(r)}_i$, i.e., $X^r_n=\{y^{(r)}_1, y^{(r)}_2, \ldots, y^{(r)}_n\}$ and $y^{(r)}_1<y^{(r)}_2<\cdots<y^{(r)}_n$. Since $D_X\not\cong D^r_n$ and $x_1\equiv y^{(r)}_1\pmod{2}$, there exists some $i_r\in [n]\backslash\{1\}$ such that $x_{i_r}\not\equiv y^{(r)}_{i_r}\pmod{2}$ for each $r=1,2$. We further assume $s_r$ to be the smallest such integer, i.e., for each $r=1,2$, $x_k\equiv y^{(r)}_k\pmod{2}$ for all $k<s_r$ and $x_{s_r}\not\equiv y^{(r)}_{s_r}\pmod{2}$. 
\noindent\par If $s_1\ge 3$ and $s_1\equiv 1\pmod{4}$ or $s_1\equiv 3\pmod{4}$, then $x_{s_1-2}\equiv x_{s_1-1}\equiv x_{s_1}\pmod{2}$. If $s_1\ge 3$ and $s_1\equiv 0\pmod{4}$ or $s_1\equiv 2\pmod{4}$, then $x_{s_1-2}\equiv x_{s_1}\not\equiv x_{s_1-1}\pmod{2}$. Now, consider the case in which $s_1=2$. Since $s_1=2$, $x_2\not\equiv y^{(1)}_2\equiv 1 \pmod{2}$ and so $x_2\equiv 2\equiv y^{(2)}_2\pmod{2}$. It follows from $x_2\equiv 2\equiv y^{(2)}_2\pmod{2}$ that $s_2\ge 3$. Now, if $s_2\equiv 0\pmod{4}$ or $s_2\equiv 2\pmod{4}$, then $x_{s_2-2}\equiv x_{s_2-1}\equiv x_{s_2}\pmod{2}$. If $s_2\equiv 1\pmod{4}$ or $s_2\equiv 3\pmod{4}$, then $x_{s_2-2}\equiv x_{s_2}\not\equiv x_{s_2-1}\pmod{2}$.
\qed

\begin{thm}\label{thmA3.7}
Let $D$ be a bipartite tournament of order $n$. There exists a tournament completion of $D$ with exactly one augmented $3$-dicycle if and only if $D$ is acyclic and $D\not\cong D^i_n$ for all $i=1,2$.
\end{thm}
\noindent\textit{Proof}: $(\Rightarrow)$ By Lemma \ref{lemA3.1}, $D$ is acyclic. By Theorem \ref{thmA2.2}, $D\cong D_X$ for some $X\subset \mathbb{Z}^+$, where $X=\{x_1, x_2, \ldots, x_n\}$ and $x_1<x_2<\cdots<x_n$. Let $T$ be a tournament completion of $D$ with exactly one augmented $3$-dicycle. If the augmented $3$-dicycle is an augmented $(2,1)$-dicycle, then by Theorem \ref{thmA3.4}, there exists some $l\in [n-2]$ such that $x_l\equiv x_{l+2}\not\equiv x_{l+1}\pmod{2}$. Then, $D\not\cong D^i_n$ for all $i=1,2$ by Lemma \ref{lemA3.6}.
\noindent\par Next, consider the case in which the augmented $3$-dicycle in $T$ is an augmented $(3,0)$-dicycle, say $u_1 u_2 u_3 u_1$, i.e., $T$ has no augmented $(2,1)$-dicycle. Suppose $D_X\cong D^i_n$ for some $i=1,2$, for a contradiction. Let $\{u_1, u_2, u_3\}$ correspond to a set $\{x,y,z\}\subset X$, where $x<y<z$. Since $u_1,u_2, u_3$ belong to a common partite set, we have $x\equiv y\equiv z\pmod{2}$. By symmetry, it suffices to consider either $(u_1, u_2, u_3)=(x,y,z)$ or $(u_1, u_2, u_3)=(x,z,y)$. In the former case, there exists some vertex $v\in N^+_D(u_1)\cap N^-_D(u_3)$ by Lemma \ref{lemA3.5}. Then $u_1 v u_3 u_1$ is an augmented $(2,1)$-dicycle in $T$, a contradiction. In the latter case, there exists some vertex $v\in N^+_D(u_1)\cap N^-_D(u_3)$ or $w\in N^+_D(u_3)\cap N^-_D(u_2)$ by Lemma \ref{lemA3.5}. Then $u_1 v u_3 u_1$ or $u_3 w u_2 u_3$ is an augmented $(2,1)$-dicycle in $T$, a contradiction.
\noindent\par $(\Leftarrow)$ Since $D$ is acyclic, by Theorem \ref{thmA2.2}, there exists some $X=\{x_1, x_2, \ldots, x_n\}$ $\subset \mathbb{Z}^+$ such that $D\cong D_X$ and $x_1<x_2<\cdots<x_n$. By the case assumption that $D\not\cong D^i_n$ for all $i=1,2$, there exists some $l\in [n-2]$ such that $x_l\equiv x_{l+1}\equiv x_{l+2}\pmod{2}$ or $x_l\equiv x_{l+2}\not\equiv x_{l+1}\pmod{2}$ by Lemma \ref{lemA3.6}. Furthermore, $x_l\equiv j\pmod{2}$ for exactly one of $j=1,2$. Define a tournament completion $T$ of $D$ as follows: for any $x,y\in X$, (\ref{eqA3.3}) holds. By (\ref{eqA3.3}) and Lemma \ref{lemA3.3}, $x_l x_{l+1} x_{l+2} x_l$ is the only dicycle in $T$. Since $D$ is acyclic, it is the only augmented $3$-dicycle.
\qed

\section{No augmented $4$-dicycles}
In this section, we study the bipartite tournaments for which there exists some tournament completion with no augmented $4$-dicycles, i.e., $t=0$ and $k=4$ in Problem \ref{probA}. Lemmas \ref{lemA4.1} and \ref{lemA4.2} underline an unified approach to tackle the cases $\mathcal{K}=\{(2,2)\}$, $\mathcal{K}=\{(3,1)\}$, and $\mathcal{K}=\{(3,1), (2,2)\}$. Let $\mathcal{A}=\{(u,v)\mid \{u, v\}\subseteq V_i \text{ for some }i=1,2\}$. The set $\mathcal{A}$ refers to the set of arcs that may be added to the bipartite tournament $D$ to obtain a tournament completion. In each case, we define later the appropriate sets $\mathcal{S}$ and $\mathcal{I}$ to invoke these lemmas and obtain the required characterization.

\begin{lem}\label{lemA4.1}
Let $D$ be a bipartite tournament. Let $\mathcal{A}=\{(u,v)\mid \{u, v\}\subseteq V_i \text{ for some }i=1,2\}$. Let $\mathcal{S}\subseteq \mathcal{A}\times \mathcal{A}$ and $\mathcal{I}\subseteq \mathcal{A}$ be sets satisfying the following properties:
\\(A1) $((x, y), (u, v))\in \mathcal{S}$ if and only if $((v, u), (y, x))\in \mathcal{S}$.
\\(A2) For each $(u_1,v_1)\in \mathcal{A}$, $(u_1,v_1)\in \mathcal{I}$ if and only if there exist some integer $k\ge 2$ and some vertices $\{u_j, v_j\mid j\in [k]\backslash\{1\}\}\subseteq V(D)$ such that $((u_i, v_i), (u_{i+1}, v_{i+1}))\in \mathcal{S}$ for all $i\in [k]$, and $(u_{k+1}, v_{k+1})=(v_l, u_l)$ for some $l\in [k-1]$.
\\If at most one of $(u, v)$ and $(v, u)$ is in the set $\mathcal{I}$ for any two vertices $u$ and $v$ in a common partite set of $D$, then there exists a tournament completion $T$ of $D$ such that $\mathcal{C}_\mathcal{S}(T)=\emptyset$, where $\mathcal{C}_\mathcal{S}(T)=\{\{x'y', v'u'\}\subseteq A(T)\mid ((x',y'),(u',v'))\in\mathcal{S}\}$.
\end{lem}
\noindent\textit{Proof}: Let $T$ be the tournament completion of $D$ with the minimum $|\mathcal{C}_\mathcal{S}(T)|$. Suppose $\mathcal{C}_\mathcal{S}(T)\neq\emptyset$ for a contradiction. Then there exist some $(x,y), (u,v)\in\mathcal{A}$ such that $\{xy, vu\}\in \mathcal{C}_\mathcal{S}(T)$, that is, $((x,y),(u,v))\in\mathcal{S}$. Let $R=\{ v^* u^*\in A(T)\mid \exists k\in\mathbb{Z}^+, ((u_i, v_i), (u_{i+1}, v_{i+1}))\in \mathcal{S} \text{ for all }i\in [k]$, where $\{(u_j, v_j)\mid j\in [k+1]\}\subseteq \mathcal{A}$ such that $u=u_1, v=v_1, u^*=u_{k+1}$, and $v^*=v_{k+1}\}\cup\{vu\}$.
\noindent\par By the hypothesis, at most one of $(u,v)$ and $(v,u)$ belongs to $\mathcal{I}$, so $(u,v)\not\in \mathcal{I}$ or $(v,u)\not\in\mathcal{I}$. First, consider the case that $(u,v)\not\in \mathcal{I}$. Define another tournament completion $T'$ of $D$ as follows: For any $u', v'\in V(T')$,
\begin{align*}
&u'v'\in A(T')\text{ if }v'u'\in R,\text{ and}\\
&v'u'\in A(T')\text{ if }v'u'\in A(T)\backslash R,
\end{align*}
i.e., we reverse all the arcs in $R$ to obtain $T'$ from $T$.
\noindent\par We claim that $|\mathcal{C}_\mathcal{S}(T')|<|\mathcal{C}_\mathcal{S}(T)|$; and this contradicts the minimality of $T$. Note that $\{xy, vu\}\in \mathcal{C}_\mathcal{S}(T)\backslash\mathcal{C}_\mathcal{S}(T')$. To show $\mathcal{C}_\mathcal{S}(T')\subset \mathcal{C}_\mathcal{S}(T)$, suppose that there exist some $(a,b), (c,d)\in\mathcal{A}$ such that $\{ab, dc\}\in\mathcal{C}_\mathcal{S}(T')\backslash\mathcal{C}_\mathcal{S}(T)$ and $((a,b),(c,d))\in\mathcal{S}$ for a contradiction. Then, exactly one of the following holds: (I) $ab\in A(T')\backslash A(T)$ and $dc\in A(T')\cap A(T)$; or (II) $dc\in A(T')\backslash A(T)$ and $ab\in A(T')\cap A(T)$; or (III) $ab, dc\in A(T')\backslash A(T)$.
\noindent\par If (I) holds, then $ba\in R$ and $((a,b),(c,d))\in\mathcal{S}$ imply that $dc\in R$. Thus, $cd\in A(T'),$ a contradiction. Note that $((a,b),(c,d))\in\mathcal{S}\Rightarrow ((d,c),(b,a))\in\mathcal{S}$ by (A1). If (II) holds, then $cd\in R$ and $((d,c),(b,a))\in\mathcal{S}$ imply that $ab\in R$. Thus, $ba\in A(T')$, a contradiction. Therefore, (III) holds and $\{ba, cd\}\subseteq R$. By the definition of $R$, there exist some positive integer $k$ (resp. $l$) and some set $\{(a_i, b_i)\mid i\in [k+1]\}\subseteq \mathcal{A}$ (resp. $\{(d_j, c_j)\mid j\in [l+1]\}\subseteq \mathcal{A}$) such that $((a_i, b_i), (a_{i+1}, b_{i+1}))\in\mathcal{S}$ for all $i\in [k]$ (resp. $((d_j, c_j), (d_{j+1}, c_{j+1}))\in\mathcal{S}$ for all $j\in [l]$), where $u=a_1=d_1$, $v=b_1=c_1$, $a=a_{k+1}$, $b=b_{k+1}$, $d=d_{l+1}$, and $c=c_{l+1}$. By (A1), $((c_{j+1}, d_{j+1}), (c_j, d_j))\in\mathcal{S}$ for all $j\in [l]$. For clarity in what follows, we denote $((p,q), (r,s))\in \mathcal{S}$ by $(p,q)\overset{\mathcal{S}}\rightsquigarrow (r,s)$. Consequently, $(u,v)=(a_1, b_1)\overset{\mathcal{S}}\rightsquigarrow (a_2, b_2)\overset{\mathcal{S}}\rightsquigarrow\cdots \overset{\mathcal{S}}\rightsquigarrow(a_{k+1}, b_{k+1})=(a,b)\overset{\mathcal{S}}\rightsquigarrow (c,d)=(c_{l+1}, d_{l+1})\overset{\mathcal{S}}\rightsquigarrow (c_l, d_l)\overset{\mathcal{S}}\rightsquigarrow\cdots \overset{\mathcal{S}}\rightsquigarrow (c_1, d_1)=(v,u)$. So, $(u,v)\in \mathcal{I}$, a contradiction. Hence, the claim follows.
\noindent\par Next, consider the case that $(v,u)\not\in \mathcal{I}$. Necessarily, $(y,x)\not\in \mathcal{I}$. Otherwise, the fact that $((x,y),(u,v))\in\mathcal{S}$, or equivalently by (A1), $((v, u), (y, x))\in \mathcal{S}$ implies that $(v,u)\in \mathcal{I}$ by (A2), a contradiction. Thus, an argument similar to the case in which $(u,v)\not\in \mathcal{I}$ can be applied on $(y,x)$, in place of $(u,v)$.
\qed

\begin{lem}\label{lemA4.2}
Let $D$ be a bipartite tournament and let $\mathcal{A}$, $\mathcal{S}$ and $\mathcal{I}$ be as given in Lemma \ref{lemA4.1}. If $(u,v)\in\mathcal{I}$ and $T$ is a tournament completion of $D$ such that $uv\in A(T)$ for some vertices $u, v\in V(D)$, then $|\mathcal{C}_\mathcal{S}(T)|>0$, where $\mathcal{C}_\mathcal{S}(T)=\{\{x'y', v'u'\}\subseteq A(T)\mid ((x',y'),(u',v'))\in\mathcal{S}\}$.
\end{lem}
\noindent\textit{Proof}: Denote $u=u_1$ and $v=v_1$. Since $(u_1,v_1)\in\mathcal{I}$, and by (A2), we may let $j$ be the maximum integer satisfying $u_r v_r\in A(T)$ for all $r\in [j]$, i.e., $\{u_j v_j, v_{j+1} u_{j+1}\}\subseteq A(T)$ and $((u_j, v_j), (u_{j+1}, v_{j+1}))\in\mathcal{S}$. Thus, $\{u_j v_j, v_{j+1} u_{j+1}\}\in\mathcal{C}_\mathcal{S}(T)$.
\qed

\indent\par
\begin{figure}[h]
\begin{subfigure}{.2\textwidth}
\begin{center}
\tikzstyle{every node}=[circle, draw, fill=black!100,
                       inner sep=0pt, minimum width=5pt]
\begin{tikzpicture}[thick,scale=0.6]%
\draw(-4,0)node[label={[yshift=-0.1cm]270:{$u_1$}}](1){};
\draw(-2,0)node[label={[yshift=-0.1cm]270:{$u_2$}}](2){};
\draw(-4,2)node[label={[yshift=0.1cm]90:{$v_1$}}](3){};
\draw(-2,2)node[label={[yshift=0.1cm]90:{$v_2$}}](4){};

\draw[->, line width=0.3mm, >=latex, shorten <= 0.1cm, shorten >= 0.1cm](1)--(2);
\draw[->, line width=0.3mm, >=latex, shorten <= 0.1cm, shorten >= 0.1cm](2)--(3);
\draw[->, line width=0.3mm, >=latex, shorten <= 0.1cm, shorten >= 0.1cm](3)--(4);
\draw[->, line width=0.3mm, >=latex, shorten <= 0.1cm, shorten >= 0.1cm](4)--(1);
\end{tikzpicture}
{\caption{$Y_1\cong C_4$}}
\end{center}
\end{subfigure}\hfill%
\begin{subfigure}{.2\textwidth}
\begin{center}
\tikzstyle{every node}=[circle, draw, fill=black!100,
                       inner sep=0pt, minimum width=5pt]
\begin{tikzpicture}[thick,scale=0.6]%
\draw(-4,0)node[label={[yshift=-0.1cm]270:{$u_1$}}](1){};
\draw(-2,0)node[label={[yshift=-0.1cm]270:{$u_2$}}](2){};
\draw(-4,2)node[label={[yshift=0.1cm]90:{$v_1$}}](3){};
\draw(-2,2)node[label={[yshift=0.1cm]90:{$v_2$}}](4){};

\draw[->, line width=0.3mm, >=latex, shorten <= 0.1cm, shorten >= 0.1cm](1)--(2);
\draw[->, line width=0.3mm, >=latex, shorten <= 0.1cm, shorten >= 0.1cm](1)--(4);
\draw[->, line width=0.3mm, >=latex, shorten <= 0.1cm, shorten >= 0.1cm](3)--(2);
\draw[->, line width=0.3mm, >=latex, shorten <= 0.1cm, shorten >= 0.1cm](3)--(4);
\end{tikzpicture}
{\caption{$Y_2$}}
\end{center}
\end{subfigure}\hfill%
\begin{subfigure}{.2\textwidth}
\begin{center}
\tikzstyle{every node}=[circle, draw, fill=black!100,
                       inner sep=0pt, minimum width=5pt]
\begin{tikzpicture}[thick,scale=0.6]%
\draw(-4,0)node[label={[yshift=-0.1cm]270:{$u_1$}}](1){};
\draw(-2,0)node[label={[yshift=-0.1cm]270:{$u_2$}}](2){};
\draw(-4,2)node[label={[yshift=0.1cm]90:{$v_1$}}](3){};
\draw(-2,2)node[label={[yshift=0.1cm]90:{$v_2$}}](4){};

\draw[->, line width=0.3mm, >=latex, shorten <= 0.1cm, shorten >= 0.1cm](2)--(1);
\draw[->, line width=0.3mm, >=latex, shorten <= 0.1cm, shorten >= 0.1cm](1)--(4);
\draw[->, line width=0.3mm, >=latex, shorten <= 0.1cm, shorten >= 0.1cm](3)--(2);
\draw[->, line width=0.3mm, >=latex, shorten <= 0.1cm, shorten >= 0.1cm](3)--(4);
\end{tikzpicture}
{\caption{$Y_3$}}
\end{center}
\end{subfigure}\hfill%
\begin{subfigure}{.2\textwidth}
\begin{center}
\tikzstyle{every node}=[circle, draw, fill=black!100,
                       inner sep=0pt, minimum width=5pt]
\begin{tikzpicture}[thick,scale=0.6]%
\draw(-4,0)node[label={[yshift=-0.1cm]270:{$u_1$}}](1){};
\draw(-2,0)node[label={[yshift=-0.1cm]270:{$u_2$}}](2){};
\draw(-4,2)node[label={[yshift=0.1cm]90:{$v_1$}}](3){};
\draw(-2,2)node[label={[yshift=0.1cm]90:{$v_2$}}](4){};

\draw[->, line width=0.3mm, >=latex, shorten <= 0.1cm, shorten >= 0.1cm](2)--(1);
\draw[->, line width=0.3mm, >=latex, shorten <= 0.1cm, shorten >= 0.1cm](4)--(1);
\draw[->, line width=0.3mm, >=latex, shorten <= 0.1cm, shorten >= 0.1cm](3)--(2);
\draw[->, line width=0.3mm, >=latex, shorten <= 0.1cm, shorten >= 0.1cm](3)--(4);
\end{tikzpicture}
{\caption{$Y_4$}}
\end{center}
\end{subfigure}\hfill%
{\caption{All non-isomorphic orientations of $K(2,2)$.}\label{figA4.5}}
\end{figure}
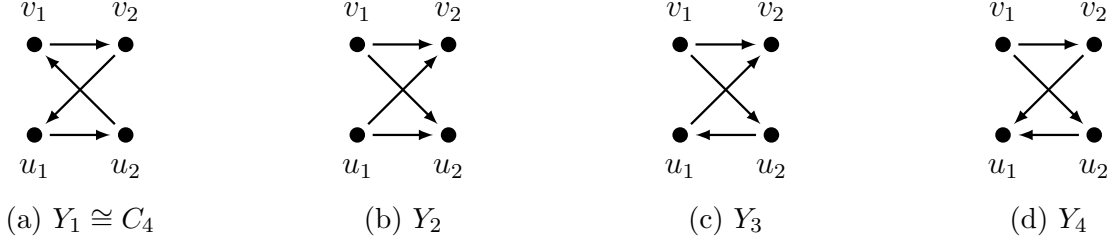

 We start by investigating the case $t=0$, $k=4$, and $\mathcal{K}=\{(2,2)\}$ of Problem \ref{probA}. There are four non-isomorphic orientations of the bipartite tournament with exactly two vertices in each partite set, namely $K(2,2)$; and we denote them by $Y_1, Y_2, Y_3$, and $Y_4$ (see Figure \ref{figA4.5}).

The verification of the next lemma is straightforward and thus, left to the reader.

\begin{lem}\label{lemA4.3}
Let $D$ be an orientation of $K(2,2)$. If there exists a tournament completion $T$ of $D$ with some augmented $(2,2)$-dicycles, then $D\in \{Y_3, Y_4\}$. Furthermore, $u_1 v_1, v_2 u_2\in A(T)$ if $D\cong Y_3$; and $u_1 v_1\in A(T)$ if $D\cong Y_4$, where the vertices are as labelled in Figure \ref{figA4.5}.
\end{lem}

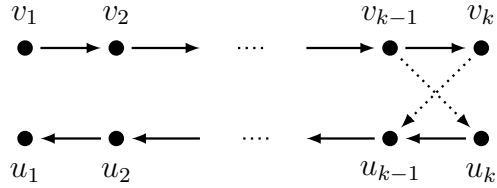
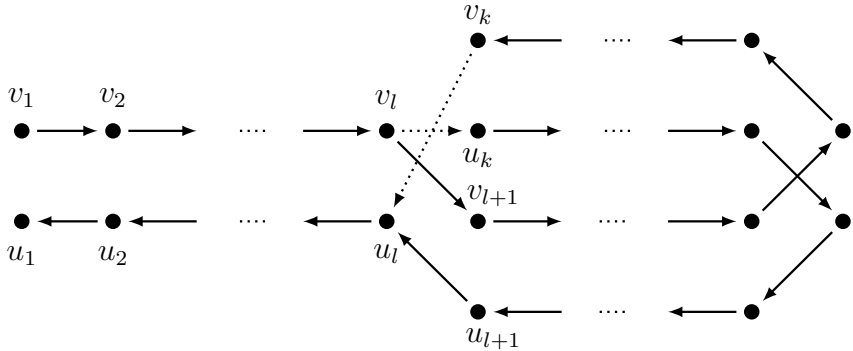
\begin{figure}[h]
\begin{center}
\begin{subfigure}{.8\textwidth}
\begin{center}
\tikzstyle{every node}=[circle, draw, fill=black!100,
                       inner sep=0pt, minimum width=5pt]
\begin{tikzpicture}[thick,scale=0.6]%
\draw(-4,4)node[label={[yshift=0.1cm]90:{$v_1$}}](v1){};
\draw(-4,2)node[label={[yshift=-0.1cm]270:{$u_1$}}](u1){};
\draw(-2,4)node[label={[yshift=0.1cm]90:{$v_2$}}](v2){};
\draw(-2,2)node[label={[yshift=-0.1cm]270:{$u_2$}}](u2){};

\draw(4,4)node[label={[yshift=-0.1cm]90:{$v_{k-1}$}}](vk-1){};
\draw(4,2)node[label={[yshift=0.1cm]270:{$u_{k-1}$}}](uk-1){};
\draw(6,4)node[label={[yshift=0.1cm]90:{$v_k$}}](vk){};
\draw(6,2)node[label={[yshift=-0.1cm]270:{$u_k$}}](uk){};

\draw[->, line width=0.3mm, >=latex, shorten <= 0.1cm, shorten >= 0.1cm](u2)--(u1);
\draw[->, line width=0.3mm, >=latex, shorten <= 0.1cm, shorten >= 0.1cm](v1)--(v2);
\draw[->, line width=0.3mm, >=latex, shorten <= 0.1cm, shorten >= 0.1cm](0,2)--(u2);
\draw[->, line width=0.3mm, >=latex, shorten <= 0.1cm, shorten >= 0.1cm](v2)--(0,4);

\draw[dotted, line width=0.3mm, >=latex, shorten <= 0.1cm, shorten >= 0.1cm](0.5,4)--(1.5,4);
\draw[dotted, line width=0.3mm, >=latex, shorten <= 0.1cm, shorten >= 0.1cm](1.5,2)--(0.5,2);

\draw[->, line width=0.3mm, >=latex, shorten <= 0.1cm, shorten >= 0.1cm](uk-1)--(2,2);
\draw[->, line width=0.3mm, >=latex, shorten <= 0.1cm, shorten >= 0.1cm](2,4)--(vk-1);

\draw[->, line width=0.3mm, >=latex, shorten <= 0.1cm, shorten >= 0.1cm](uk)--(uk-1);
\draw[->, line width=0.3mm, >=latex, shorten <= 0.1cm, shorten >= 0.1cm](vk-1)--(vk);
\draw[dotted, ->, line width=0.3mm, >=latex, shorten <= 0.1cm, shorten >= 0.1cm](vk)--(uk-1);
\draw[dotted, ->, line width=0.3mm, >=latex, shorten <= 0.1cm, shorten >= 0.1cm](vk-1)--(uk);
\end{tikzpicture}
{\caption{$(u_1, v_1)$ is $d$-inconsistent and $l=k-1$, i.e., $(u_{k+1}, v_{k+1})=(v_{k-1}, u_{k-1})$.}}
\end{center}
\end{subfigure}
\vfill%
\begin{subfigure}{.8\textwidth}
\begin{center}
\tikzstyle{every node}=[circle, draw, fill=black!100,
                       inner sep=0pt, minimum width=5pt]
\begin{tikzpicture}[thick,scale=0.6]%
\draw(-4,4)node[label={[yshift=0.1cm]90:{$v_1$}}](v1){};
\draw(-4,2)node[label={[yshift=-0.1cm]270:{$u_1$}}](u1){};
\draw(-2,4)node[label={[yshift=0.1cm]90:{$v_2$}}](v2){};
\draw(-2,2)node[label={[yshift=-0.1cm]270:{$u_2$}}](u2){};

\draw(4,4)node[label={[yshift=0.1cm]90:{$v_l$}}](vl){};
\draw(4,2)node[label={[yshift=-0.1cm]270:{$u_l$}}](ul){};

\draw(6,6)node[label={[yshift=0cm, xshift=0cm]90:{$v_k$}}](vk){};
\draw(6,4)node[label={[yshift=0cm, xshift=0cm]270:{$u_k$}}](uk){};

\draw(12,6)node[label={[yshift=0cm]90:{}}](5){};
\draw(12,4)node[label={[yshift=0cm]270:{}}](6){};

\draw(6,2)node[label={[yshift=-0.15cm, xshift=0.2cm]90:{$v_{l+1}$}}](vl+1){};
\draw(6,0)node[label={[yshift=0.15cm, xshift=0.2cm]270:{$u_{l+1}$}}](ul+1){};

\draw(12,2)node[label={[yshift=-0.25cm, xshift=-0.2cm]90:{}}](1){};
\draw(12,0)node[label={[yshift=0.2cm, xshift=-0.2cm]270:{}}](2){};

\draw(14,4)node[label={[xshift=0.1cm]0:{}}](3){};
\draw(14,2)node[label={[xshift=0.1cm]0:{}}](4){};

\draw[->, line width=0.3mm, >=latex, shorten <= 0.1cm, shorten >= 0.1cm](u2)--(u1);
\draw[->, line width=0.3mm, >=latex, shorten <= 0.1cm, shorten >= 0.1cm](v1)--(v2);
\draw[->, line width=0.3mm, >=latex, shorten <= 0.1cm, shorten >= 0.1cm](0,2)--(u2);
\draw[->, line width=0.3mm, >=latex, shorten <= 0.1cm, shorten >= 0.1cm](v2)--(0,4);

\draw[dotted, line width=0.3mm, >=latex, shorten <= 0.1cm, shorten >= 0.1cm](1.5,2)--(0.5,2);
\draw[dotted, line width=0.3mm, >=latex, shorten <= 0.1cm, shorten >= 0.1cm](1.5,4)--(0.5,4);

\draw[->, line width=0.3mm, >=latex, shorten <= 0.1cm, shorten >= 0.1cm](ul)--(2,2);
\draw[->, line width=0.3mm, >=latex, shorten <= 0.1cm, shorten >= 0.1cm](2,4)--(vl);

\draw[dotted, ->, line width=0.3mm, >=latex, shorten <= 0.1cm, shorten >= 0.1cm](vl)--(uk);
\draw[dotted, ->, line width=0.3mm, >=latex, shorten <= 0.1cm, shorten >= 0.1cm](vk)--(ul);

\draw[dotted, line width=0.3mm, >=latex, shorten <= 0.1cm, shorten >= 0.1cm](9.5,4)--(8.5,4);
\draw[dotted, line width=0.3mm, >=latex, shorten <= 0.1cm, shorten >= 0.1cm](9.5,6)--(8.5,6);

\draw[->, line width=0.3mm, >=latex, shorten <= 0.1cm, shorten >= 0.1cm](uk)--(8,4);
\draw[->, line width=0.3mm, >=latex, shorten <= 0.1cm, shorten >= 0.1cm](8,6)--(vk);

\draw[->, line width=0.3mm, >=latex, shorten <= 0.1cm, shorten >= 0.1cm](5)--(10,6);
\draw[->, line width=0.3mm, >=latex, shorten <= 0.1cm, shorten >= 0.1cm](10,4)--(6);

\draw[->, line width=0.3mm, >=latex, shorten <= 0.1cm, shorten >= 0.1cm](ul+1)--(ul);
\draw[->, line width=0.3mm, >=latex, shorten <= 0.1cm, shorten >= 0.1cm](vl)--(vl+1);
\draw[->, line width=0.3mm, >=latex, shorten <= 0.1cm, shorten >= 0.1cm](8,0)--(ul+1);
\draw[->, line width=0.3mm, >=latex, shorten <= 0.1cm, shorten >= 0.1cm](vl+1)--(8,2);

\draw[dotted, line width=0.3mm, >=latex, shorten <= 0.1cm, shorten >= 0.1cm](8.5,2)--(9.5,2);
\draw[dotted, line width=0.3mm, >=latex, shorten <= 0.1cm, shorten >= 0.1cm](8.5,0)--(9.5,0);

\draw[->, line width=0.3mm, >=latex, shorten <= 0.1cm, shorten >= 0.1cm](10,2)--(1);
\draw[->, line width=0.3mm, >=latex, shorten <= 0.1cm, shorten >= 0.1cm](2)--(10,0);

\draw[->, line width=0.3mm, >=latex, shorten <= 0.1cm, shorten >= 0.1cm](1)--(3);
\draw[->, line width=0.3mm, >=latex, shorten <= 0.1cm, shorten >= 0.1cm](4)--(2);

\draw[->, line width=0.3mm, >=latex, shorten <= 0.1cm, shorten >= 0.1cm](3)--(5);
\draw[->, line width=0.3mm, >=latex, shorten <= 0.1cm, shorten >= 0.1cm](6)--(4);
\end{tikzpicture}
{\caption{$(u_1, v_1)$ is $d$-inconsistent and $(u_{k+1}, v_{k+1})=(v_l, u_l)$ for some $l\in [k-3]$.}}
\end{center}
\end{subfigure}%
\end{center}
{\caption{$(u_1, v_1)$ is $d$-inconsistent.}\label{figA4.6}}
\end{figure}

\begin{defn}\label{defnA4.4}
Let $D$ be a bipartite tournament. Let $\mathcal{A}=\{(u,v)\mid \{u, v\}\subseteq V_i$ for some $i=1,2\}$ and suppose $\{(u_i, v_i)\mid i\in [r+1]\}\subseteq \mathcal{A}$ for some $r\in\mathbb{Z}^+$.
\\(a) We say that the ordered pair $(u_1, v_1)$ \textit{d-specifies} the ordered pair $(u_2, v_2)$ if $u_1\leftarrow u_2$ and $v_1\rightarrow v_2$. Here, the ``$d$" in ``$d$-specifies" indicates that $u_1$ and $v_1$ belong to the same partite set, while $u_2$ and $v_2$ belong to another partite set, which is \textit{different} from the first one.
\\(b) We say that the ordered pair $(u_1, v_1)$ is \textit{$d$-inconsistent} if there exists some integer $2\le k\le r$ such that $(u_i, v_i)$ $d$-specifies $(u_{i+1}, v_{i+1})$ for all $i\in [k]$, and $(u_{k+1}, v_{k+1})=(v_l, u_l)$ for some $l\in [k-1]$ (see Figure \ref{figA4.6}(a) and (b) when $l=k-1$ and $l\in [k-3]$, respectively).
\end{defn}

\begin{rmk}\label{rmkA4.5}
Clearly, if $(u_1, v_1)$ $d$-specifies $(u_2, v_2)$ and $T$ is a tournament completion of $D$ with $u_1 v_1, v_2 u_2\in A(T)$, then $u_1 v_1 v_2 u_2 u_1$ is an augmented $(2,2)$-dicycle in $T$.
\end{rmk}

\begin{rmk} \label{rmkA4.6}
Let us make some remarks on Definition \ref{defnA4.4}(b).
~\\(a) We emphasize that the definition does not require the vertices $u_i$'s and $v_j$'s to be all pairwise distinct, i.e., it is possible to have $u_i=u_j$, or $v_i=v_j$, or $u_i=v_j$ for some distinct $i,j\in [k]$.
\\(b) For the case in which $l\in [k-3]$ (and $k\ge 4$), $k-l$ is at least three and necessarily odd since $D$ is a bipartite tournament.
\\(c) For each $i\in [k]$, either $D\langle \{u_i, v_i, u_{i+1}, v_{i+1}\}\rangle \cong Y_3$ or $D\langle \{u_i, v_i, u_{i+1}, v_{i+1}\}\rangle \cong Y_4$.
\end{rmk}

\begin{ppn}\label{ppnA4.7}
Let $D$ be a bipartite tournament. Then, there exists a tournament completion of $D$ with no augmented $(2,2)$-dicycles if and only if for any two vertices $u$ and $v$ in a common partite set of $D$, at most one of $(u, v)$ and $(v, u)$ is $d$-inconsistent.
\end{ppn}
\noindent\textit{Proof}: Let $\mathcal{A}=\{(u,v)\mid \{u, v\}\subseteq V_i \text{ for some }i=1,2\}$. Define the sets $\mathcal{S}\subseteq \mathcal{A}\times \mathcal{A}$ and $\mathcal{I}\subseteq \mathcal{A}$ such that $\mathcal{S}=\{((x,y),(u,v))\in\mathcal{A}\times \mathcal{A}\mid (x,y)$ $d$-specifies $(u,v)\}$ and $\mathcal{I}=\{(u,v)\in\mathcal{A}\mid (u,v)$ is $d$-inconsistent$\}$. It is straightforward to verify that $\mathcal{S}$ and $\mathcal{I}$ satisfy the properties (A1) and (A2) in Lemma \ref{lemA4.1}.
\noindent\par $(\Rightarrow)$ Suppose that $(u,v)$ and $(v,u)$ are $d$-inconsistent for some vertices $u, v\in V(D)$. Let $T$ be a tournament completion of $D$. By Lemma \ref{lemA4.2}, regardless of whether $uv\in A(T)$ or $vu\in A(T)$, there exist some arcs $x'y', v'u'\in A(T)$ such that $(x',y')$ $d$-specifies $(u',v')$. By Remark \ref{rmkA4.5}, $T$ contains some augmented $(2,2)$-dicycle.
\noindent\par $(\Leftarrow)$  By the assumption and Lemma \ref{lemA4.1}, there exists a tournament completion $T$ of $D$ such that $\{\{x'y', v'u'\}\subseteq A(T)\mid ((x',y'),(u',v'))\in\mathcal{S}\}=\emptyset$, i.e., $T$ contains no augmented $(2,2)$-dicycles.
\qed

\begin{cor}
Let $D$ be a bipartite tournament containing no $Y_3$. Then there exists a tournament completion of $D$ with no augmented $(2,2)$-dicycles if and only if $D$ does not contain $F$ (given in Figure \ref{figA1.1}).
\end{cor}
\noindent\textit{Proof}: $(\Rightarrow)$ Suppose $D$ contains $F$ as labelled in Figure \ref{figA1.1}. It is easy to check that $(u_1, u_2)$ and $(u_2, u_1)$ are both $d$-inconsistent. By Proposition \ref{ppnA4.7}, every tournament completion of $D$ has some augmented $(2,2)$-dicycle.
\noindent\par $(\Leftarrow)$ Suppose $T$ is a tournament completion of $D$ that contains some augmented $(2,2)$-dicycle. By Proposition \ref{ppnA4.7}, there exist two vertices $u$ and $v$ in a common partite set of $D$ such that both $(u, v)$ and $(v, u)$ are $d$-inconsistent. Since $D$ contains no $Y_3$, it follows from Lemma \ref{lemA4.3} and Remark \ref{rmkA4.6}(c) that $D$ contains $F$.
\qed

\indent\par
\begin{figure}[h]
\begin{subfigure}{.55\textwidth}
\begin{center}
\tikzstyle{every node}=[circle, draw, fill=black!100,
                       inner sep=0pt, minimum width=5pt]
\begin{tikzpicture}[thick,scale=0.6]%
\draw(-4,2)node[label={[xshift=-0.1cm]180:{$u_1$}}](u1){};
\draw(-4,4)node[label={[xshift=-0.1cm]180:{$v_1=v_2$}}](v1){};
\draw(-4,6)node[label={[xshift=-0.1cm]180:{$u_2$}}](u2){};
\draw(-2,4)node[label={[xshift=0cm]0:{$w$}}](w){};

\draw[->, line width=0.3mm, >=latex, shorten <= 0.1cm, shorten >= 0.1cm](w)--(u1);
\draw[->, line width=0.3mm, >=latex, shorten <= 0.1cm, shorten >= 0.1cm](u2)--(w);

\end{tikzpicture}
{\caption{$v_1=v_2$ and $d_D(u_2, u_1)=2$}}
\end{center}
\end{subfigure}%
\begin{subfigure}{.4\textwidth}
\begin{center}
\tikzstyle{every node}=[circle, draw, fill=black!100,
                       inner sep=0pt, minimum width=5pt]
\begin{tikzpicture}[thick,scale=0.6]%
\draw(-4,2)node[label={[xshift=-0.1cm]180:{$v_2$}}](v2){};
\draw(-4,4)node[label={[xshift=-0.1cm]180:{$u_1=u_2$}}](u1){};
\draw(-4,6)node[label={[xshift=-0.1cm]180:{$v_1$}}](v1){};
\draw(-2,4)node[label={[xshift=0cm]0:{$w$}}](w){};

\draw[->, line width=0.3mm, >=latex, shorten <= 0.1cm, shorten >= 0.1cm](v1)--(w);
\draw[->, line width=0.3mm, >=latex, shorten <= 0.1cm, shorten >= 0.1cm](w)--(v2);

\end{tikzpicture}
{\caption{$u_1=u_2$ and $d_D(v_1, v_2)=2$}}
\end{center}
\end{subfigure}\hfill%
{\caption{$(u_1, v_1)$ $c$-specifies $(u_2, v_2)$.}\label{figA4.7}}
\end{figure}
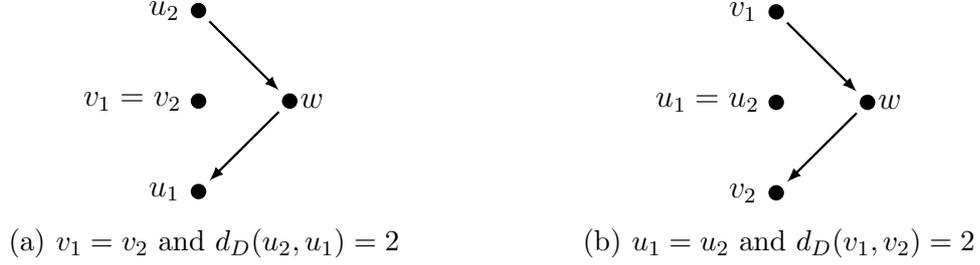

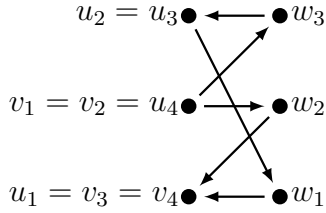
\begin{figure}[h]
\begin{center}
\tikzstyle{every node}=[circle, draw, fill=black!100,
                       inner sep=0pt, minimum width=5pt]
\begin{tikzpicture}[thick,scale=0.6]%
\draw(2,2)node[label={[yshift=0cm]180:{$u_1=v_3=v_4$}}](1){};
\draw(2,4)node[label={[yshift=0cm]180:{$v_1=v_2=u_4$}}](2){};
\draw(2,6)node[label={[yshift=0cm]180:{$u_2=u_3$}}](3){};

\draw(4,2)node[label={[yshift=0cm]0:{$w_1$}}](4){};
\draw(4,4)node[label={[yshift=0cm]0:{$w_2$}}](5){};
\draw(4,6)node[label={[yshift=0cm]0:{$w_3$}}](6){};

\draw[->, line width=0.3mm, >=latex, shorten <= 0.1cm, shorten >= 0.1cm](4)--(1);
\draw[->, line width=0.3mm, >=latex, shorten <= 0.1cm, shorten >= 0.1cm](5)--(1);

\draw[->, line width=0.3mm, >=latex, shorten <= 0.1cm, shorten >= 0.1cm](2)--(5);
\draw[->, line width=0.3mm, >=latex, shorten <= 0.1cm, shorten >= 0.1cm](2)--(6);

\draw[->, line width=0.3mm, >=latex, shorten <= 0.1cm, shorten >= 0.1cm](6)--(3);
\draw[->, line width=0.3mm, >=latex, shorten <= 0.1cm, shorten >= 0.1cm](3)--(4);
\end{tikzpicture}
\captionsetup{justification=centering}
{\caption{Subdigraph of a bipartite tournament in which $(u_1, v_1)$ is $c$-inconsistent.\\Note that $(u_i, v_i)$ $c$-specifies $(u_{i+1}, v_{i+1})$ for all $i\in [4]$ and $(u_4,v_4)=(v_1, u_1)$.}\label{figA4.8}}
\end{center}
\end{figure}

 Our next aim is to characterize the bipartite tournaments for which there exists a tournament completion with no augmented $(3,1)$-dicycles. To this end, we introduce the following notions, which are required for defining the sets $\mathcal{S}$ and $\mathcal{I}$ later in Proposition \ref{ppnA4.11}.
\begin{defn}\label{defnA4.9}
Let $D$ be a bipartite tournament. Let $\mathcal{A}=\{(u,v)\mid \{u, v\}\subseteq V_i$ for some $i=1,2\}$ and suppose $\{(u_i, v_i)\mid i\in [r+1]\}\subseteq \mathcal{A}$ for some $r\in\mathbb{Z}^+$.
\\(a) We say that the ordered pair $(u_1, v_1)$ \textit{$c$-specifies} the ordered pair $(u_2, v_2)$ if $v_1=v_2$ and $d_D(u_2, u_1)=2$ or $u_1=u_2$ and $d_D(v_1, v_2)=2$ (see Figure \ref{figA4.7}). Here, the ``$c$" in ``$c$-specifies" indicates that the vertices $u_1, v_1, u_2$, and $v_2$ are in a \textit{common} partite set.
\\(b) We say that the ordered pair $(u_1, v_1)$ is \textit{$c$-inconsistent} if there exists some integer $2\le k\le r$ such that $(u_i, v_i)$ $c$-specifies $(u_{i+1}, v_{i+1})$ for all $i\in [k]$, and $(u_{k+1}, v_{k+1})=(v_l, u_l)$ for some $l\in [k-1]$ (see Figure \ref{figA4.8} for an example).
\end{defn}

\begin{rmk}\label{rmkA4.10} Suppose $\{u_1, v_1, u_2, v_2\}\subseteq V_j$ for some $j=1,2$. Clearly, if $(u_1, v_1)$ $c$-specifies $(u_2, v_2)$ and $T$ is a tournament completion of $D$ with $u_1 v_1, v_2 u_2\in A(T)$, then either $u_1 v_1 u_2 w u_1$ or $v_2 u_1 v_1 w v_2$ is an augmented $(3,1)$-dicycle in $T$ for some $w\in V_{3-j}$.
\end{rmk}

\begin{ppn}\label{ppnA4.11}
Let $D$ be a bipartite tournament. There exists a tournament completion of $D$ with no augmented $(3,1)$-dicycles if and only if for any two vertices $u$ and $v$ in a common partite set of $D$, at most one of $(u, v)$ and $(v, u)$ is $c$-inconsistent.
\end{ppn}
\noindent\textit{Proof}: Let $\mathcal{A}=\{(u,v)\mid \{u, v\}\subseteq V_i \text{ for some }i=1,2\}$. Define the sets $\mathcal{S}\subseteq \mathcal{A}\times \mathcal{A}$ and $\mathcal{I}\subseteq \mathcal{A}$ such that $\mathcal{S}=\{((x,y),(u,v))\in\mathcal{A}\times \mathcal{A}\mid (x,y)$ $c$-specifies $(u,v)\}$ and $\mathcal{I}=\{(u,v)\in\mathcal{A}\mid (u,v)$ is $c$-inconsistent$\}$. It is straightforward to verify that $\mathcal{S}$ and $\mathcal{I}$ satisfy the properties (A1) and (A2) in Lemma \ref{lemA4.1}.
\noindent\par $(\Rightarrow)$ Suppose that $(u,v)$ and $(v,u)$ are $c$-inconsistent for some vertices $u, v\in V(D)$. Let $T$ be a tournament completion of $D$. By Lemma \ref{lemA4.2}, regardless of whether $uv\in A(T)$ or $vu\in A(T)$, there exist some arcs $x'y', v'u'\in A(T)$ such that $(x',y')$ $c$-specifies $(u',v')$. By Remark \ref{rmkA4.10}, $T$ contains some augmented $(3,1)$-dicycle.
\noindent\par $(\Leftarrow)$  By the assumption and Lemma \ref{lemA4.1}, there exists a tournament completion $T$ of $D$ such that $\{\{x'y', v'u'\}\subseteq A(T)\mid ((x',y'),(u',v'))\in\mathcal{S}\}=\emptyset$, i.e., $T$ contains no augmented $(3,1)$-dicycles.
\qed

\indent\par The last theorem of this section, namely Theorem \ref{thmA4.13}, generalizes Propositions \ref{ppnA4.7} and \ref{ppnA4.11} to characterize the bipartite tournaments for which there exists a tournament completion with no augmented $\{(3,1), (2,2)\}$-dicycles. To achieve this, we generalize Definitions \ref{defnA4.4} and \ref{defnA4.9} as follows.

\begin{defn}\label{defnA4.12}
Let $D$ be a bipartite tournament. Let $\mathcal{A}=\{(u,v)\mid \{u, v\}\subseteq V_i$ for some $i=1,2\}$ and suppose $\{(u_i, v_i)\mid i\in [r+1]\}\subseteq \mathcal{A}$ for some $r\in\mathbb{Z}^+$.
\\(a) We say that the ordered pair $(u_1, v_1)$ \textit{specifies} the ordered pair $(u_2, v_2)$ if $(u_1, v_1)$ $c$-specifies or $d$-specifies $(u_2, v_2)$.
\\(b) We say that the ordered pair $(u_1, v_1)$ is \textit{inconsistent} if there exists some integer $k\ge 2$ such that $(u_i, v_i)$ specifies $(u_{i+1}, v_{i+1})$ for all $i\in [k]$, and $(u_{k+1}, v_{k+1})=(v_l, u_l)$ for some $l\in [k-1]$.
\end{defn}

\begin{thm}\label{thmA4.13}
Let $D$ be a bipartite tournament and $\mathcal{K}=\{(3,1), (2,2)\}$. There exists a tournament completion of $D$ with no augmented $\mathcal{K}$-dicycles if and only if for any two vertices $u$ and $v$ in a common partite set of $D$, at most one of $(u, v)$ and $(v, u)$ is inconsistent.
\end{thm}
\noindent\textit{Proof}: Let $\mathcal{A}=\{(u,v)\mid \{u, v\}\subseteq V_i \text{ for some }i=1,2\}$. Define the sets $\mathcal{S}\subseteq \mathcal{A}\times \mathcal{A}$ and $\mathcal{I}\subseteq \mathcal{A}$ such that $\mathcal{S}=\{((x,y),(u,v))\in\mathcal{A}\times \mathcal{A}\mid (x,y)$ specifies $(u,v)\}$ and $\mathcal{I}=\{(u,v)\in\mathcal{A}\mid (u,v)$ is inconsistent$\}$. It is straightforward to verify that $\mathcal{S}$ and $\mathcal{I}$ satisfy the properties (A1) and (A2) in Lemma \ref{lemA4.1}.
\noindent\par $(\Rightarrow)$ Suppose that $(u,v)$ and $(v,u)$ are inconsistent. Let $T$ be a tournament completion of $D$. By Lemma \ref{lemA4.2}, regardless of whether $uv\in A(T)$ or $vu\in A(T)$, there exist some arcs $x'y', v'u'\in A(T)$ such that $(x',y')$ specifies $(u',v')$. By Remarks \ref{rmkA4.5} and \ref{rmkA4.10}, $T$ contains some augmented $\mathcal{K}$-dicycle.
\noindent\par $(\Leftarrow)$  By the assumption and Lemma \ref{lemA4.1}, there exists a tournament completion $T$ of $D$ such that $\{\{x'y', v'u'\}\subseteq A(T)\mid ((x',y'),(u',v'))\in\mathcal{S}\}=\emptyset$, i.e., $T$ contains no augmented $\mathcal{K}$-dicycles.
\qed

\section{Conclusion}
In this paper, we examine Problem \ref{probA} with a focus on tournament completions with exactly one augmented $3$-dicycle (i.e., $t=1$ and $k=3$ as denoted in Problem \ref{probA}) and with no $4$-dicycles (i.e., $t=0$ and $k=4$). The case $k=3$ depends mainly on a characterization of acyclic bipartite tournaments, namely Theorem \ref{thmA2.2}, and a characterization of acyclic bipartite tournaments not isomorphic to $D^1_n$ and $D^2_n$, namely Lemma \ref{lemA3.6}. The approach for the case of $k=4$ is established on Lemmas \ref{lemA4.1} and \ref{lemA4.2}, and identifying the appropriate notions (Definitions \ref{defnA4.4}, \ref{defnA4.9}, and \ref{defnA4.12}) for $\mathcal{S}$ and $\mathcal{I}$ in each subcase $\mathcal{K}=\{(2,2)\}$, $\mathcal{K}=\{(3,1)\}$, and $\mathcal{K}=\{(3,1), (2,2)\}$, respectively. We hope that this work illustrates the potential for future research in other values of $t$ and $k$. Results concerning dicycles of tournaments and bipartite tournaments are often extended to multipartite tournaments. Therefore, if we change ``bipartite tournaments" to ``multipartite tournaments" in Problem \ref{probA}, it becomes a more general problem. 

\section*{Acknowledgement}
The author would like to thank the editors and referees for their helpful comments.


\begin{thebibliography}{99}
\bibitem{BA CT}
{B. Alspach and C. Tabib,}{ A note on the number of 4-circuits in a tournament}, \textit{Ann. Discrete Math.}, 12, (1982), 13--19.
\bibitem{KSB LWB}
{K.S. Bagga and L.W. Beineke,}{ On the numbers of some subtournaments of a bipartite tournament}, \textit{Ann. New York Acad. Sci.}, 555, (1989), 21--29.
\bibitem{JB GG}
{J. Bang-Jensen and G. Gutin,}\textit{ Digraphs: Theory, Algorithms and Applications}, Springer Verlag, London, (2009).
\bibitem{JB JH XZ}
{J. Bang-Jensen, J. Huang and X. Zhu,}{ Completing orientations of partially oriented graphs}, \textit{J. Graph Theory}, 87, (2017), 285--304.
\bibitem{LWB CHCL}
{L.W. Beineke and C.H.C. Little,}{ Cycles in bipartite tournaments}, \textit{J. Combin. Theory Ser. B}, 32, (1982), 140--145.
\bibitem{SD PG SG SS}
{S. Das, P. Ghosh, S. Ghosh and S. Sen,}{ Oriented bipartite graphs and the Goldbach graph}, \textit{Discrete Math.}, 344, (2021), Article 112497.
\bibitem{FH LM}
{F. Harary and L. Moser,}{ The theory of round robin tournaments}, \textit{Amer. Math. Monthly}, 73, (1966), 231--246.
\bibitem{KH and JH}
{K. Hsu and J. Huang,}{ Obstructions for local tournament orientation completions}, \textit{Discrete Math.}, 346, (2023), Article 113220.
\bibitem{NL AM}
{N. Linial and A. Morgenstern,}{ On the number of 4-cycles in a tournament}, \textit{J. Graph Theory}, 83, (2015), 266--276.
\bibitem{JWM2}
{J.W. Moon,}\textit{ Topics on tournaments}, {Holt, Rinehart and Winston, New York}, (1968).
\bibitem{LV}
{L. Volkmann,}{ Cycles in multipartite tournaments: results and problems}, \textit{Discrete Math.}, 245, (2002), 19--53.
\end{thebibliography}
\end{document}